\documentclass[11pt]{amsart}

\usepackage{amssymb}
\usepackage{listings}

\usepackage{tikz}
\usepackage[ruled,vlined,linesnumbered,norelsize]{algorithm2e}
\SetEndCharOfAlgoLine{}
\SetKwProg{Function}{Function}{}{}
\SetKwInOut{Input}{Input}
\SetKwInOut{Output}{Output}

\usepackage{url}

\newcommand{\Z}{\mathbb{Z}}
\newcommand{\Q}{\mathbb{Q}}
\newcommand{\R}{\mathbb{R}}
\newcommand{\K}{\mathbf{K}}

\newcommand{\mat}[4]{\left(\begin{array}{cc} {#1} & {#2} \\ {#3} & {#4}\end{array}\right)}
\newcommand{\GFq}[1]{\mathbb{F}_{#1}}
\newcommand{\EE}{\mathcal{E}}
\newcommand{\Es}{\mathcal{E}^*}
\newcommand{\As}{A^*}
\newcommand{\Bs}{B^*}
\newcommand{\ComputeUMod}{{\sc ComputeUMod}}
\newcommand{\ComputeAMod}{{\sc ComputeAMod}}

\newcommand{\powersums}{\mathcal{P}}
\newcommand{\sfM}{{\sf M}}
\newcommand{\qdiff}[1]{{#1}'}

\newcommand{\Fell}[1]{F_{#1}}

\newcommand{\GLtwo}{\mathrm{GL}_2}

\newtheorem{lemma}{Lemma}[section]
\newtheorem{theorem}[lemma]{Theorem}
\newtheorem{proposition}[lemma]{Proposition}
\newtheorem{corollary}[lemma]{Corollary}

\begin{document}

\title{Using Fricke modular polynomials \\ to
compute isogenies}
\author{François Morain}
\address{
    LIX - Laboratoire d'informatique de l'École polytechnique 
    \emph{and}
    GRACE - Inria Saclay--Île-de-France
}
\email{morain@lix.polytechnique.fr} 

\date{\today}

\maketitle

\begin{abstract}
Let $\EE$ be an elliptic curve over a field $\K$ and $\ell$ a prime.
There exists an elliptic curve $\Es$ related to $\EE$ by an
isogeny of degree $\ell$ only if $\Phi_\ell^t(X, j(\EE)) = 0$, where
$\Phi_\ell^t(X, Y)$ is the traditional modular polynomial.
Moreover, $\Phi_\ell^t$
gives the coefficients of $\Es$, together with parameters needed to
build the isogeny explicitly. Since $\Phi_\ell^t$
has very large coefficients, many families with smaller
coefficients can be used instead, as described by Elkies, Atkin and others.
In this work, we concentrate on the computation of the family of
modular polynomials introduced by Fricke and more recently used
by Charlap, Coley and Robbins. In some cases, the resulting
polynomials are small, which justifies the interest of this study.
We review and adapt the known algorithms to perform the computations
of these polynomials. After describing the use of series computations, we
investigate fast algorithms using floating point numbers based on fast
numerical evaluation of Eisenstein series. We also explain how to use
isogeny volcanoes as an alternative. The
last part is concerned with finding explicit formulas for computing
the coefficients of $\Es$. To this we add tables of numerical examples.
\end{abstract}

%%%%% S
\section{Introduction}

Computing isogenies is the central ingredient of the
Schoof-Elkies-Atkin (SEA) algorithm
that computes the cardinality of elliptic curves over finite fields of
large characteristic~\cite{Schoof95,Atkin92b,Elkies98} and also
\cite{BlSeSm99}.
More recently, it has found its way in post-quantum
cryptography~\cite{ChLaGo09,JaFe11,CaLaMaPaRe18,FeKoLePeWe20} among
others, as well as the cryptosystems~\cite{Couveignes06,RoSt06,FeKiSm18}.

Let $\K$ be a field of characteristic different from $2$ and $3$.
A (separable) isogeny between two elliptic curves $\EE/\K: Y^2 = X^3 +
A X + B$ and
$\Es/\K: Y^2 = X^3 + \As X + \Bs$ is a group morphism that is a rational
map of degree $N$ (the cardinality of its kernel assumed to be cyclic).
There are two ways to handle these isogenies. When the degree $N$
is small, formulas for $\As$, $\Bs$ and the kernel polynomial can be
precomputed. For large $N$,
one of the key ingredients is modular polynomials, the second one
finding rational expressions for $\As$ and $\Bs$ from $A$, $B$ and the
modular polynomial. Note there is a purely
algebraic approach using triangular sets~\cite[\S 7]{PoSc13} (see also
\cite{NoYaYo20}).

There are many families of modular polynomials that can be used, with
different properties. Very generally, a modular polynomial is some
bivariate polynomial $\Phi(X, J)$ where $J$ corresponds to the
$j$-invariant of the elliptic curve $\EE$, and $X$ stands for some
modular function on  $\Gamma_0(N)$.
The prototype is $j(\Es)$ (see below for
more precise statements) that
yields traditional modular polynomials. Alternative choices for $X$
exist. They all lead to polynomials of (conjectured) height $O(\psi(N)
\log N)$ but with small constants. Here $\psi(N) = N \prod_{p \mid N}
(1+1/p)$ is the cardinality of $\Gamma_0(N) \backslash \Gamma$.

In \cite{Fricke22}, Fricke computes a resolvant polynomial $U_N$ of
degree $\psi(N)$ for $\wp(N z)$ where $\wp$ is the Weierstrass
$\wp$-function of $\EE$ and $N$ small. It turns out to be the same
polynomial as used in Elkies's work and also by \cite{ChCoRo91}
(without notice). 
The latter authors complete this with two polynomials
$V_N$, $W_N$ having the property (among others) that $\As$
(resp. $\Bs$) is a root of $V_N(X, A, B)$ (resp. $W_N(X, A,
B)$). The polynomial $U_N$ is the modular polynomial associated
with a form of weight 2 for $\Gamma_0(N)$ related to the Eisenstein
series $E_2$; $V_N$ (resp. $W_N$) is the modular polynomial for
$E_4(N \tau)$ (resp. $E_6(N \tau)$).

The aim of this work is to describe the properties of the
polynomials $(U_N, V_N, W_N)$ and the relevant algorithms to
compute them, extending methods already used in the traditional
cases of bivariate modular polynomials based on modular
functions. These methods apply to modular forms of even weight and
yield trivariate polynomials in $(E_4, E_6, \Delta)$, see below for
the rationale. Though our primary interest is in prime
values $N = \ell$, the theory and practice for the general case of an
integer $N$ follow the same paths. In Appendix~\ref{ex6}, we give some
details on an example when $N = 6$.

Section~\ref{sct:formulas}) recalls results on classical functions and
modular forms, adding approaches to recognize a collection of modular
forms as polynomials in the traditional quantities $E_4$, $E_6$ and
$\Delta$. Section~\ref{sct:fast} addresses fast methods to evaluate
these functions, including an approach to the
simultaneous evaluation of several classical series.
Section~\ref{sct:Fricke} gathers properties
of Fricke polynomials;
Section~\ref{sct:computations} is devoted to
various algorithms for performing the computations developed for
classical modular polynomials.
This includes the computation of algebraic expressions for $\As$ and
$\Bs$ as rational fractions (see
\cite{NoYaYo20}). Section~\ref{sct:FrickeA} explains how to compute the
isogenous curve using partial derivatives of the polynomial $U_\ell$,
in the spirit of Atkin's work in the traditional case.
We give numerical examples and height comparisons between modular
polynomials in Section~\ref{sct:results}. An appendix contains
numerical values for our polynomials, as well as a script for checking
the results of Section~\ref{sct:FrickeA}.

\newcommand{\qt}{q_1}
\bigskip
\noindent
{\bf Notations:} Let $\tau$ be a complex number in the upper half
plane. Put $\qt = \exp(i \pi \tau)$ and $q = \qt^2$. Depending on
authors, formulas are expressed in either parameter, which sometimes is
clumsy. We write indifferently $f(\tau)$ or $f(q)$ some series. We
denote by $n^\omega$ the complexity of multiplying two $n\times n$
matrices over some field, that is used in the analysis of computer
algebra algorithms, in our case solving linear systems.

%%%%% S
\section{Modular forms}
\label{sct:formulas}

For convenience, we follow~\cite{CoSt17} and the references given
below are related to this.

%%%%% SSS
\subsection{Elements of theory}

%%%%%%%%%% SS
\subsubsection{Integer matrices}

We note:
$$\Gamma = \mathrm{Sl}_2(\Z) = \left\{
\mat{a}{b}{c}{d}, a d - b c = 1\right\}.$$
Let $n > 0$ be an integer and
$\mathcal{M}_n$ be the set of primitive integral $2\times 2$
matrices of determinant $n$. The relation $R \sim_\Gamma R'$ means
$\Gamma R = \Gamma R'$ for $R$, $R'$ in $\mathcal{M}_n$.

\begin{proposition}[Prop. 6.5.3]
The matrices
$$R_i = \mat{a_i}{b_i}{0}{d_i}, 1 \leq i \leq \psi(n)$$
with $a_i > 0$, $a_i d_i = n$, $\gcd(a_i, b_i, d_i)=1$ and $0 \leq b_i
< d_i$, form a set of representatives of
$\mathcal{M}_n$ modulo $\sim_\Gamma$.
\end{proposition}

\begin{proposition}[Prop. 6.5.3 -- cont'd]\label{prop:ellm}
Let $M = \mat{a}{b}{c}{d}$ be of determinant $N$. Then
$$M \sim_\Gamma R = \mat{A}{B}{0}{D}$$
where $A = \gcd(a, c)$ and $D = N/A$.
\end{proposition}

\medskip
\noindent
{\em Proof:} We look for $V = \mat{u}{v}{w}{x} \in
\Gamma$ such that $V M = \mat{A}{B}{0}{D}$ for $AD=N$, $\gcd(A, B,
D)=1$. This is equivalent to
\begin{eqnarray}
u a + v c &=& A, \label{eqna1} \\
w a + x c &=& 0, \label{eqna2} \\
u b + v d &=& B, \label{eqna3} \\
w b + x d &=& D.\label{eqna4}
\end{eqnarray}
A natural candidate for $A$ is $\gcd(a, c) = u a + v c$ in integers
$u$ and $v$. Note that $A \mid N$ since $ad-bc=N$. Write $a = A
a'$, $c = A c'$ with $\gcd(a', c')=1$.
Equation (\ref{eqna2}) has solutions $w = -c'$ and $x = a'$, which
is coherent with equation (\ref{eqna4}), with the choice of $D =
N/A$. We are left with $B_0 = u b + v d = B + q D$,
with $0 \leq B < D$ so that
$$\mat{1}{-q}{0}{1} \cdot V M =
\mat{1}{-q}{0}{1} \mat{A}{B + q D}{0}{D}
= \mat{A}{B}{0}{D}.$$
We finish with
$$U = V^{-1} \cdot \mat{1}{q}{0}{1} = \mat{a'}{a' q - v}{c'}{c' q + u}$$
and we are done. $\Box$

We denote by
$$\Gamma_0(N) = \left\{ \mat{a}{b}{c}{d} \in \Gamma, c \equiv 0 \bmod
N\right\},$$
and $\psi(N) = N \prod_{p \mid N} (1+1/p)$ the cardinal of
$\Gamma_0(N) \backslash \Gamma$.

The general case is treated in Corollary~6.2.11 but we need the prime
case only.
\begin{proposition}[Ex. 6.2.12]\label{prop:cosets}
Let $N = \ell$ be a prime number.
A system of representative of cosets for $\Gamma_0(\ell) \backslash
\Gamma$ is formed of the
matrices $R_c = \mat{1}{0}{c}{1}$ for $0 \leq c < \ell$, and $R_\ell =
\mat{0}{-1}{1}{\ell}$.
\end{proposition}
The following result will prove useful for computing conjugate values
of the functions that we study later on.
\begin{lemma}\label{cbar}
For $0 < c < \ell$:
$$\mat{\ell}{0}{c}{1} = \mat{\ell}{-\overline{c}}{c}{u}
\mat{1}{\overline{c}}{0}{\ell} = T_c \; R_c'$$
with $u \ell + \overline{c}{c} = 1$;
$$\mat{0}{-\ell}{1}{\ell} = \mat{0}{-1}{1}{1} \;
\mat{1}{0}{0}{\ell} = T_\ell \; R_\ell'.$$
\end{lemma}

%%%%%%%%%% SS
\subsubsection{Modular forms}

Let $f$ be a modular form for $\Gamma$ of integer weight $2k$. By
construction
$$f(M \tau) = (c \tau + d)^{2k} f(\tau)$$
for any matrix $M = \mat{a}{b}{c}{d}$ in $\Gamma$.

If $M \in \GLtwo(\R)$, define
$$f |_{2k} M(\tau) = \det(M)^{k} (c \tau + d)^{-2k} f(M \tau).$$

%%%%% SSS
\subsection{Eisenstein series}

The classical Eisenstein series\footnote{Ramanujan used $L = P = E_2$,
$M = Q = E_4$, $N = R = E_6$.} we consider are
$$E_2(q) = 1 -24 \sum_{n=1}^\infty \delta_1(n) q^n,$$
$$E_4(q) = 1 + 240 \sum_{n=1}^\infty \delta_3(n) q^n,$$
$$E_6(q) = 1 - 504 \sum_{n=1}^\infty \delta_5(n)q^n,$$
where $\delta_r(n)$ denotes the sum of the $r$-th powers of the
divisors of $n$. Other series $E_{2k}$ can be defined for even $k > 3$.
The series $E_{2k}$ is a modular form of weight $2k$ for $k
> 1$.

Also of interest is the discriminant $\Delta$:
$$\Delta(q) = (E_4(q)^3-E_6(q)^2)/1728 = \eta(q)^{24}.$$
Dedekind's function is $\eta(q) = q^{1/24} \prod_{n=1}^\infty
(1-q^n)$.

Finally, the modular invariant is
$$j(q) = \frac{E_4(q)^3}{\Delta(q)} = \frac{1}{q} + 744 + \cdots.$$

The series $E_2$ is not a modular form since (see \cite{Rankin77} or
Corollary 5.2.17):
\begin{theorem}\label{E2tsf}
For all matrices
$\left(\begin{array}{cc} a & b \\ c & d \\\end{array}\right)$
in $\Gamma$, one has
\begin{equation}\label{eqE2}
E_2((a\tau+b)/(c\tau+d)) = (c\tau+d)^2 E_2(\tau) + \frac{6 c}{\pi i} (c \tau
+ d).
\end{equation}
\end{theorem}

We can build a modular form easily as follows. Let $N$ be an integer and
let $\Fell{N}$ denote the {\em multiplier} $E_2(\tau)-N
E_2(N\tau)$. From \cite{Ramanujan14} and \cite{Berndt89}, we get
\begin{proposition}
The function $\Fell{N}$ is a modular form of weight 2 and
trivial multiplier system for $\Gamma_0(N)$.
\end{proposition}

\noindent
{\em Proof:} Write, for $ad-bc=1$ and $N \mid c$, the value
$$E_2\left(N \frac{a\tau+b}{c\tau+d}\right)
= E_2\left(\frac{a (N\tau)+N b}{(c/N)(N\tau)+d}\right)
= ((c/N)N\tau+d)^2 E_2(N\tau) - \frac{6 c}{N \pi i}
((c/N) N\tau + d)$$
which leads to
$$N E_2\left(N \frac{a\tau+b}{c\tau+d}\right)
= N (c\tau+d)^2 E_2(N\tau) - \frac{6 c}{\pi i}
(c\tau + d).$$
Subtracting $E_2((a\tau+b)/(c\tau+d))$, we see that
$$\Fell{N}((a\tau+b)/(c\tau+d)) = (c\tau+d)^2
\Fell{N}(\tau). \Box$$

Some identities are known for small values of $N$, for instance
\cite{KaKo03} for $N \in \{2, 4\}$; \cite{Berndt89} for $N=3$ and $11$;
\cite[Thm 6.2]{BoBo91}, 
\cite[Thm 3.7]{BeChSoSo00} for $N \in \{5, 7\}$. A very nice
relation is~\cite[Thm 6.3]{BeChSoSo00}
$$\Fell{7}(q) = 6 \left(\sum_{m, n = -\infty}^\infty
q^{m^2+mn+2n^2}\right)^2.$$

%%%%%%%%%%%%%%% SSS
\subsection{Formulas}

Reference is Proposition~2.4.1.
When $f(q) = \sum_{n\geq n_0} a_n q^n$, we introduce the operator
\begin{equation}\label{operator}
\qdiff{f}(q) = \frac{1}{2 i \pi} \, \frac{df}{d\tau} = q \frac{d f}{d
q} = \sum_{n \geq n_0} n a_n q^n.
\end{equation}
Several identities are classical:
\begin{equation}\label{jdelta}
\Delta = \frac{E_4^3-E_6^2}{1728}, \quad
j = \frac{E_4^3}{\Delta}, \quad j - 1728 = \frac{E_6^2}{\Delta},
\end{equation}
\begin{equation}\label{formj}
\frac{\qdiff{j}}{j} = - \frac{E_6}{E_4}, \quad \frac{\qdiff{j}}{j-1728} = -
\frac{E_4^2}{E_6}, \quad \qdiff{j} = - \frac{E_4^2 E_6}{\Delta},
\quad \frac{\qdiff{\Delta}}{\Delta} = E_2,
\end{equation}
to which we add the Ramanujan differential system:
\begin{equation}\label{diff46}
{3 \qdiff{E_4}} = {E_4} E_2 -
{E_6}, \quad {2 \qdiff{E_6}} = {E_6} E_2 - {E_4^2},
\quad 12 \qdiff{E_2} = E_2^2 - E_4.
\end{equation}

%%%%%%%%%%%%%%% SSS
\subsection{Expressing modular forms as polynomials}
\label{sss:express}

We extend a remark already done in \cite{ChCoRo91} that uses the
following result from \cite[\S 5.6.2]{Serre73}.
\begin{theorem}
A modular form $f$ of weight $2k$ with integer coefficients can be
expressed as a polynomial with integer coefficients in $E_4$, $\Delta$
if $k$ is even and $E_4$, $\Delta$, $E_6$ otherwise. The number of
coefficients in this polynomial is approximately $k/6$.
\end{theorem}
We can apply this result to higher index Eisenstein series. For instance
\begin{equation}\label{Eeq}
E_8 = E_4^2, E_{10} = E_4 E_6.
\end{equation}

\begin{lemma}\label{lem1}
Let $k > 1$. Consider the equation
\begin{equation}\label{eq4612}
2 i_4 + 3 i_6 + 6 i_{12} = k
\end{equation}
where $i_4$, $i_{12}$ are positive integers and $i_6 \in \{0, 1\}$.
Write $k = 2 k_0 + \epsilon$ for $\epsilon \in \{0, 1\}$ and $m = k_0
- \epsilon$.
All solutions $(i_6, i_4, i_{12})$ to equation~(\ref{eq4612}) are
$$(\epsilon, m-3j, j) \text{ for } 0\leq j \leq j_{\max} := \lfloor
m/3\rfloor.$$
\end{lemma}

\medskip
\noindent
{\em Proof:} If $k$ is even, we write $k = 2m$, which forces $i_6 = 0$.
We rewrite
(\ref{eq4612}) as $i_4 + 3 i_{12} = m$ and the result follows.
If $k = 2m+1$, we need $i_6 = 1$ and
(\ref{eq4612}) becomes $i_4 + 3 i_{12} = m-1$, which concludes the
proof. $\Box$

\medskip
\noindent
{\em Proof of the theorem:} the expression we are looking for is
$$f = \sum_{i_4, i_6, i_{12}}
c_{i_4, i_6, i_{12}} E_6^{i_6} E_4^{i_4} \Delta^{i_{12}}$$
for all positive indices satisfying
$4 i_4 + 6 i_6 + 12 i_{12} = 2 k$ and $i_6 \in \{0, 1\}$.
which is in fact equivalent to (\ref{eq4612}) and we apply the
Lemma. $\Box$

For $w = 2k$, write
$$P_{w, j} = E_6^{\epsilon} E_4^{m-3j} \Delta^{j} = q^{j} +
\sum_{n > j} g_{j, n} q^n$$
where the coefficients $g_{j, n}$ are integers and do not depend on
$f$. We may precompute all $P_{w, j}$'s for $0 \leq j \leq j_{\max} =
\lfloor m/3\rfloor$ using Algorithm~\ref{algoPwj}.
\LinesNotNumbered
\begin{algorithm}[hbt]
\caption{Evaluating all $P_{w, j}$'s. \label{algoPwj}}
\SetKwProg{Fn}{Function}{}{}
\Fn{EvaluateAllPwj($E_4, E_6^\epsilon, \Delta, w, j_{\max}$)}{
\Input{$E_4$, $E_6$, $\Delta$; $w = 2k$, $j_{\max}$}
\Output{$(P_{w, j})$ for $0 \leq j \leq j_{\max}$}

0. compute $m$\;

1. [compute $P_{w, j} = E_6^\epsilon \cdot \Delta^j$ for all $j$]

1.1. $P_{w, 0} \leftarrow E_6^\epsilon$\;

1.2. \For{$j \leftarrow 1$ \KwTo $j_{\max}$}{
  $P_{w, j} \leftarrow \Delta \cdot P_{w, j-1}$\;
}

2. $Q[0] \leftarrow 1$; $Q[1] \leftarrow E_4$; $Q[2] \leftarrow
E_4^2$; $Q[3] \leftarrow E_4 \cdot Q[2]$\;

3. $m' = m-3 j_{\max}$; $S \leftarrow Q[m']$\;

4. \For{$j \leftarrow j_{\max-1}$ \KwTo $0$}{
  $S \leftarrow Q[3] \cdot S$\; \tcp{$S = E_4^{m-3j}$}
  $P_{w, j} \leftarrow P_{w, j} \cdot S$\;
}
\Return{$\{P_{w, j}\}$}\;
}
\end{algorithm}
The cost of this algorithm is $3 j_{\max} + 2$ multiplications of 
series. We could improve on this using
squarings in Step 1.2. If this is a one-time computation, Step 4 may not
update and store the $P_{w, j}$'s.

The next two results are crucial for our forthcoming computations.
\begin{proposition}\label{prop1}
Let $f = \sum_{n \geq 0} f_n q^n$ be a modular form of weight
$w = 2k$. With the notations of Lemma~\ref{lem1}, there exist numbers
$c_j$ such that
$$f = \sum_{j=0}^{j_{\max}} c_j P_{w, j}$$
where the $c_j$'s are solution of a triangular linear system.
\end{proposition}

\medskip
\noindent
{\em Proof:} Write
$$\sum_{n \geq 0} f_n q^n = \sum_{j=0}^{j_{\max}}
c_{j} \left(q^{j} + \sum_{n > j} g_{j, n} q^n\right),$$
the system giving the $c_j$'s is triangular. Solving the system takes
$O(j_{\max}^2) = O(w^2)$ operations with a very small constant. $\Box$

From this, we deduce:
\begin{corollary}
When the $f_n$'s are integers, so are the $c_{j}$'s.
\end{corollary}

\begin{corollary}\label{order}
To express $f$ as a polynomial, we need all series developped at order
$j_{\max} \approx w/12$. 
\end{corollary}

To prepare for the computation of modular equations of
Section~\ref{sct:computations}, we need to express a collection of
modular forms of weight $2k r$ for $1\leq r\leq \psi(N)$ ($\ell$ is an
odd prime $>3$) as polynomials in $E_4$, $E_6$ and $\Delta$. The
rationale is to share evaluations of products $P_{\epsilon, x, y} =
E_6^\epsilon E_4^x \Delta^y$. It follows from Corollary~\ref{order}
that the forms must be computed with order close to $w \psi(N) / 12 =
k \psi(N)/6$.

\begin{figure}[hbt]
\begin{center}
\begin{tikzpicture}[>=latex,scale=0.25]
\node at (1, 0) {$\bullet$};
\node at (2, 0) {$\bullet$};
\node at (3, 0) {$\bullet$};
\node at (0, 1) {$\bullet$};
\node at (4, 0) {$\bullet$};
\node at (1, 1) {$\bullet$};
\node at (5, 0) {$\bullet$};
\node at (2, 1) {$\bullet$};
\node at (6, 0) {$\bullet$};
\node at (3, 1) {$\bullet$};
\node at (0, 2) {$\bullet$};
\draw[step=1cm] (0, 0) grid (6, 2);
\node at (3, -2) {$\ell = 5$, $k = 2$};
\end{tikzpicture}
\hspace*{3cm}
\begin{tikzpicture}[>=latex,scale=0.25]
\node at (3, 0) {$\bullet$};
\node at (0, 1) {$\bullet$};
\node at (6, 0) {$\bullet$};
\node at (3, 1) {$\bullet$};
\node at (0, 2) {$\bullet$};
\node at (9, 0) {$\bullet$};
\node at (6, 1) {$\bullet$};
\node at (3, 2) {$\bullet$};
\node at (0, 3) {$\bullet$};
\node at (12, 0) {$\bullet$};
\node at (9, 1) {$\bullet$};
\node at (6, 2) {$\bullet$};
\node at (3, 3) {$\bullet$};
\node at (0, 4) {$\bullet$};
\draw[step=1cm] (0, 0) grid (12, 4);
\node at (6, -2) {$\ell=7$, $k=3$, $\epsilon=0$};
\end{tikzpicture}
\hspace*{1cm}
\begin{tikzpicture}[>=latex,scale=0.25]
\node at (0, 0) {$\bullet$};
\node at (3, 0) {$\bullet$};
\node at (0, 1) {$\bullet$};
\node at (6, 0) {$\bullet$};
\node at (3, 1) {$\bullet$};
\node at (0, 2) {$\bullet$};
\node at (9, 0) {$\bullet$};
\node at (6, 1) {$\bullet$};
\node at (3, 2) {$\bullet$};
\node at (0, 3) {$\bullet$};
\draw[step=1cm] (0, 0) grid (9, 3);
\node at (9/2, -2) {$\ell=7$, $k=3$, $\epsilon=1$};
\end{tikzpicture}
\end{center}
\caption{Examples of grids and sets $\mathcal{I}_\epsilon$. \label{grids}}
\end{figure}
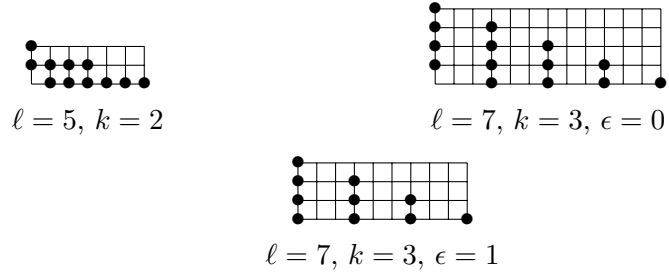

Using notations from Lemma~\ref{lem1} for $1\leq r\leq \psi(N)$, write
$k r = 2 k_r+\epsilon_r$, $m_r = k_r - \epsilon_r$. Each $m_r$ gives
rise to a set of indices $(m_r - 3 j_r, j_r)$ for $0 \leq j_r \leq m_r/3$. We
can gather all these points in the plane to give rise to interesting
patterns, as indicated in Figure~\ref{grids}. We define two sets of
indices $\mathcal{I}_\epsilon$ corresponding to pairs $(m_r - 3 j_r,
j_r)$ with $\epsilon_r = \epsilon \in \{0, 1\}$. When $k$ is even,
$\mathcal{I}_1$ is empty.

\begin{lemma}
The sets $\mathcal{I}_\epsilon$ satisfy the following properties.

i) $\mathcal{I}_0$ is included in the grid
$(0, 0) \times (m_{\max, 0}, j_{\max, 0})$ with
$m_{\max, 0} = \lfloor k \psi(N)/2\rfloor$, $j_{\max, 0} = \lfloor k
\psi(N)/6\rfloor$;

ii) the distance between two consecutive abscissas in $\mathcal{I}_0$
is $k/2$ when $k$ is even and $k$ when $k$ is odd;

iii) for any abscissa $x$ of a point in $\mathcal{I}_0$, the maximal
$y$ is $\lfloor (m_{\max, 0}-x)/3\rfloor$.

\noindent
For $k$ odd, $\mathcal{I}_1$ is such that

iv) $\mathcal{I}_1 \subset (0, 0) \times (m_{\max, 1}, j_{\max, 1})$ with
$m_{\max, 1} = k (r'-3)/2$ with $r'$ the largest odd integer $\leq
\psi(N)$ and $j_{\max, 1} = \lfloor (k \ell-3)/6\rfloor$;

v) the distance between two consecutive abscissas in $\mathcal{I}_1$
is $k$;

vi) for any abscissa $x$ in $\mathcal{I}_1$, the maximal $y$ is
$\lfloor (m_{\max, 1}-x)/3\rfloor$.
\end{lemma}

\noindent
{\em Proof:}
i) the maximal value of $y$ is always reached for
$m_{\max, 0} = m_{\psi(N)} = \lfloor k \psi(N)/2\rfloor$.

ii) When $k$ is even, all $\epsilon_r$'s are $0$,
$m_r = k_r$, the minimal value is $m_1 = k/2$. Two consecutive values
have distance $m_{r+1}-m_r = k_{r+1}-k_r = k/2$.
When $k$ is odd, the smallest $m$ is $m_2 = k$. Two consecutive values
are such that $m_{2r+2}-m_{2r} = k$.

iii) A point $(x, y)$ is in $\mathcal{I}_0$ if and only if there
exists $0 \leq m \leq m_{\max, 0}$ and $0 \leq j \leq j_{\max, 0}$
such that $(x, y) = (m-3j, j)$. It follows that $j \leq (m-x)/3 \leq
(m_{\max, 0}-x)/3$.

iv) In the remainder of the proof, $k$ is odd. When $\epsilon=1$
(which cannot happen unless both $k$ and $r$ are odd)
we get $m_r = (k r - 3)/2$, and the maximal value
is $m_{\max, 1} = (k r' -3)/2$. Note that
$r' = \ell$ when $N = \ell$ is prime.

v) The minimal value is $m_1 = \max(0, (k-3)/2)$. Two consecutive values
are such that $m_{2r+3}-m_{2r+1} = k$.

vi) Proceed as in case iii. $\Box$

We write algorithms in a generic manner, and it will
work for series, floating point numbers, etc.
The auxiliary routine performing operations is given in
 Algorithm~\ref{algoPeps}. 
 We need one multiplication for each point in
 $\mathcal{I}_\epsilon$. Each power of $E_4$ is computed
 incrementally. The main function is given in
 Algorithm~\ref{algoAllP}. The simple algorithm expressing the
 coefficients of the representation of $f$ is give as
 Algorithm~\ref{algoexpress}.

\LinesNotNumbered
\begin{algorithm}[hbt]
\caption{Compute all $P_{\epsilon, x, y}$'s. \label{algoPeps}}
\SetKwProg{Fn}{Function}{}{}
\Fn{EvaluateAllPepsilon($\epsilon$, $E_4, E_6^\epsilon, \Delta$,
$m_{\max}$, $x_0$, $dx$)}{
\Input{$E_4$, $E_6^\epsilon$ and $\Delta$, $m_{\max}$, $x_0$, $dx$}
\Output{$(P_{\epsilon, x, y})$ for all points in $\mathcal{I}_\epsilon$}

0. $j_{\max} \leftarrow \lfloor m_{\max}/3\rfloor$\;

1. compute $P_{\epsilon, 0, j}$ for all $0 \leq j \leq j_{\max}$\;

3. \For{$x \leftarrow x_0$ \KwTo $m_{\max}$ by $dx$}{
 compute $P_{\epsilon, x, y} = E_6^\epsilon E_4^x \Delta^y$ for
all $0 \leq y \leq \lfloor (m_{\max}-x)/3\rfloor$\;
  }
}
\end{algorithm}

\LinesNotNumbered
\begin{algorithm}[hbt]
\caption{Evaluating all $P_{\epsilon, x, y}$'s. \label{algoAllP}}
\SetKwProg{Fn}{Function}{}{}
\Fn{EvaluateAllP($\psi$, $k$, $E_4, E_6^\epsilon, \Delta$)}{
\Input{$E_4$, $E_6^\epsilon$ and $\Delta$, $\psi = \psi(N)$, $k$}
\Output{$(P_{\epsilon, x, y})$ for all points in
$\mathcal{I}_\epsilon$ and all $\epsilon$}

1. \If{$k$ is even}{
$dx \leftarrow k/2$\;
}
\Else{
$dx \leftarrow k$\;
}

2. $P_{0, x, y} \leftarrow$ EvaluateAllPepsilon($0$, $E_4$, $E_6$,
$\Delta$, $k \psi/2$, $dx$, $dx$)\;

3. \If{$k$ is odd}{
  $P_{1, x, y} \leftarrow$ EvaluateAllPepsilon($1$, $E_4$, $E_6$, $\Delta$,
  $k (r'-3)/2$, $(k-3)/2$, $k$) where $r'$ is the largest odd integer
  $\leq \psi$\;
  }

4. \Return{$\{P_{0, x, y}\} \cup \{P_{1, x, y}\}$}.
}
\end{algorithm}

\LinesNotNumbered
\begin{algorithm}[hbt]
\caption{Express $f$ as a polynomial in the $P_{w, j}$'s. \label{algoexpress}}
\SetKwProg{Fn}{Function}{}{}
\Fn{ExpressForm($f$, $P_{w, j}$, $j_{\max}$)}{
\Input{$f$ a modular form of weight $w$}
\Output{The coefficients of the representation of $f$ as
$\sum_{j=0}^{j_{\max}} c_j P_{w, j}$}

0. $g \leftarrow f$\;
1. \For{$j \leftarrow 0$ \KwTo $j_{\max}$}{
\tcp{$g = g_j q^j + \cdots$}
 1.1 $c_j \leftarrow g_j$\;

 1.2 $g \leftarrow g - c_j P_{w, j}$\;
  }

2. \Return{$c_j$'s}.
}
\end{algorithm}

%%%%%%%%%% SS
\subsection{Modular polynomials from modular forms}
\label{sct:mpmf}

Let us begin with modular forms $f$ for $\Gamma_0(N)$ of even weight $w$.
In this section $(R_i)$ denote a list of representatives of cosets of
$\Gamma_0(N) \backslash \Gamma$.

%%%%% SSS
\subsubsection{General results}

\begin{theorem}\label{mp1}
Let $f$ be a modular form of weight $w$ for $\Gamma_0(N)$. The values
$f|_w(R_i)$ are conjugate over $\Gamma_0(N)$ and define a polynomial
$$\Phi[f](X) = \prod_{R} (X - f|_w(R_i)) = X^{\psi(N)} + C_1(f)
X^{\psi(N)-1} + \cdots + C_{\psi(N)}(f).$$

a) The coefficient $C_t(f)$ is a modular form of weight $w t$ for
$\Gamma$.

b) The polynomial $\Phi[f](X, E_4, E_6, \Delta)$ is homogeneous with
weight $w (\ell+1)$.
\end{theorem}

\medskip
\noindent
{\em Proof:} 

a) $C_t(f)$ is a symmetric function of the
$f|_w(R_i)$'s. Having a matrix of $\Gamma$ operate on the $(R_i)$'s leave
them globally invariant, so that $C_t(f)$ is invariant under $\Gamma$.

b) is a consequence of a). $\Box$

\medskip
In practice, it is customary to compute the power sums of roots of
$\Phi[f]$:
$$S_t(f) = \sum_{i=1}^{\psi(N)} (f|_w(R_i))^t$$
for $1\leq t \leq \psi(N)$. The sum $S_t(f)$ is also a modular form of weight
$w t$ for $\Gamma$, hence is expressible as a polynomial in $(E_4,
E_6, \Delta)$. The functions $S_t(f)$ are more easily computed and once
recognized as polynomials in $(E_4, E_6, \Delta)$, we can recover the
coefficients of $\Phi[f]$ using Newton's formulas using $O({\sf
M}(\psi(N)))$ operations (in a field of characteristic 0 or larger
than $\psi(N)$).

Let us turn towards a special case.
\begin{proposition}
If $f$ is a modular form of weight $w$ for $\Gamma$, then
the function $g(\tau) = f(N \tau)$ is a modular form of weight $w$ for
$\Gamma_0(N)$.
\end{proposition}

\medskip
\noindent
{\em Proof:} let $M = \mat{a}{b}{c}{d} \in \Gamma_0(N)$, in other
words $c = N c'$ for some integer $c'$. Write
$$g(M \tau) = f(N (M\tau)) = f\left(\frac{(Na)
\tau+(Nb)}{(Nc')\tau+d}\right).$$
Consider the matrix
$$U = \mat{a}{b N}{c'}{d}$$
which belongs to $\Gamma$ since $a d - b (Nc') = 1$. Using
$$M = U \cdot \mat{N}{0}{0}{1},$$
we get
$$g(M \tau) = f(U (N\tau)) = (c' (N\tau)+d)^w f(N \tau) =
(c \tau+d)^w g(\tau),$$
which proves invariance. $\Box$

Inspired by \cite[Example 6.2.15]{CoSt17}.
\begin{theorem}
Let $(R_i)$ be a system of cosets for $\Gamma_0(N) \backslash \Gamma$, the $(g
|_w R_i)$ are permuted by $\Gamma$ and we let $\Phi[g](X)$ denote the
modular polynomial
$$\Phi[g](X) = \prod_i (X - g|_w R_i).$$
If $R = \mat{a}{b}{c}{d}$ is a representative of a coset, the matrix
$\mat{aN}{bN}{c}{d}$ is equivalent to a
matrix $\mat{A}{B}{0}{D}$ with $A = \gcd(a, c)$, $AD=N$ and $0 \leq B
< D$. Moreover
$$g|_w R(\tau) = D^{-w} f\left(\frac{A \tau + B}{D}\right).$$
\end{theorem}

\medskip
\noindent
{\em Proof:} If a coset is $R = \mat{a}{b}{c}{d}$ with $d\mid N$, $d > 0$,
$c$ defined modulo $N/d$:
\begin{eqnarray*}
g |_w R (\tau) &=& (c\tau+d)^{-w} g\left(\frac{a \tau + b}{c \tau +
d}\right) \\
&=& (c\tau+d)^{-w} f\left(\frac{a N \tau + b N}{c \tau + d}\right) \\
&=& N^{-w/2} f|_w \mat{a N}{b N}{c}{d} (\tau).
\end{eqnarray*}
The last matrix is equivalent to some $R_c' = \mat{A}{B}{0}{D}$ with $A D =
N$ and $0 \leq B < D$ given by Proposition~\ref{prop:ellm}, so that
\begin{eqnarray*}
g |_w R (\tau) &=& N^{-w/2} f|_w \mat{A}{B}{0}{D} (\tau) \\
&=& N^{-w/2} \left(N^{w/2} D^{-w} f\left(\frac{A \tau +
B}{D}\right)\right) \\
&=& D^{-w} f\left(\frac{A \tau + B}{D}\right). \quad \Box\\
\end{eqnarray*}
In order to handle series with integer coefficients, we scale all
conjugates by multiplication by $N^w$.

%%%%% SSS
\subsubsection{Orders and heights}

The traditional modular polynomial $\Phi_N^t$ has height $6
\psi(N) \log N$ approximately~\cite{Cohen84,Sutherland13}. To get a
general result for modular forms, we first estimate the order of the
series needed.

For integer $t > 0$, write
$$f^t(q) = \sum_{n=0}^\infty \alpha(t, n) q^n.$$
\begin{proposition}\label{propSDf}
Let $D \mid N$. Then
$$S_{D, t}(f) = \sum_{B=0}^{D-1} f^t\left(\frac{A \tau + B}{D}\right)
 = D \sum_{n=0}^\infty \alpha(t, D n) q^{An}.$$
\end{proposition}

\medskip
\noindent
{\em Proof:} using $z = q^{A/D}$:
$$S_{D, t}(f) = \sum_{B=0}^{D-1} \sum_{n=0}^\infty
\alpha(t, n) z^n \zeta_D^{B n}
= \sum_{n=0}^\infty \alpha(t, n) z^n \sum_{B=0}^{D-1} \zeta_D^{B n}.$$
If $n$ is prime to $D$, the inner sum $\Sigma_D(n)$ is 0 
using the properties of roots of unity. If $D \mid n$, the sum is $D$.
Now, suppose that $g = \gcd(n, D)1$ with $1 < g <$, and $n' = n/g$
which is prime to $D/g$. Write $B = B' + x D'$, $0 \leq B' < D'$. We
obtain
$$\Sigma_D(n) = \sum_{x=0}^{g-1} \sum_{B'=0}^{D'-1} \zeta_{D'}^{B' n'
+ x n' D'}
= \sum_{x=0}^{g-1} \sum_{B'=0}^{D'-1} \zeta_{D'}^{B'} = 0.$$
Finally:
$$S_{D, t}(f) = D \sum_{n=0}^\infty \alpha(t, D n) q^{An}. \quad \Box$$

\begin{proposition}\label{proporder}
If $f$ is a modular form of weight $w$ for $\Gamma$. We
need to develop $f$ up to order $w N \psi(N)/12$ to compute
$\Phi[f(N\tau)]$.
\end{proposition}

\medskip
\noindent
{\em Proof:} we need to compute the expansion of the sum
$$S_t(f) = \sum_{D \mid N} (N/D)^{wt} S_{D, t}(f).$$
To get $\mathcal{O}$ terms in the sum, we need to develop the series
up to order $N \mathcal{O}$.
If we need $\mathcal{O}$ terms, we need to have $N n \geq
\mathcal{O}$.

Since $S_t(f)$ is a modular form of weight $w t$, we need $(wt)/12$
terms with a maximal $\mathcal{O} = w \psi(N)/12$, which yields the
result. $\Box$

\smallskip
Now, we turn our attention to estimating the height of $\Phi$.
For this, we need general bounds on the coefficients of modular
forms. This is \cite[Theorem 9.2.1]{CoSt17}:
\begin{theorem}
Let $f(q) = \sum_{n \geq 0} a_n q^n$ be a modular form of weight $w$
for $\Gamma$.

a) if $f$ is a cusp form, then $a(n) = O(n^{w/2})$.

b) if $f$ is not cuspidal (i.e., $a_0 \neq
0$), then there are two positive constants
$C_1$, $C_2$ s.t.
$$C_1 n^{w-1} \leq |a(n)| \leq C_2 n^{w-1}.$$
\end{theorem}
For instance, for Eisenstein series:
$$n^{w-1} \leq \delta_{w-1}(n) \leq \zeta(n-1) n^{w-1}.$$
We could use better estimates, see \cite[Chapter 9]{CoSt17}.

\smallskip
We slightly generalize the computations in \cite{ChCoRo91} to the
case of an arbitrary form.
\begin{proposition}\label{sizeME}
Let $f$ be a modular form for $\Gamma$ of weight $w$ and $K > 0$ an
integer, $C > 0$ two constants such that $|a(n)| \leq C n^{K-1}$. The
size of the sum $S_{\psi(N)}(f)$ is approximately $K \psi(N) \log (N
\psi(N))$.
\end{proposition}

\medskip
\noindent
{\em Proof:} Let us denote $\psi(N)$ by $\psi$ for short.
With $u = |q|$, we get
$$|f(q)| \leq \sum_{n=0}^\infty n^{K-1} u^n
\leq C \sum_{n=0}^\infty (n+K-1) (n+K-2) \cdots (n+1) u^n = C (K-1)!
(1-u)^{-K} := C_1 (1-u)^{-K},$$
from which $|f^m(q)| \leq C_1^m (1-u)^{-Km}$ for all $m$.
By Proposition~\ref{proporder}, we need to compute $f^\psi$ at order
$d \approx w N \psi/12$.
The largest coefficient is therefore of the order of
$$L = C_1^\psi \frac{(d+K \psi-1)!}{d!}.$$
Using Stirling's formula, we get
$$\log L \approx \psi \log C_1 + 
d \log\left(\frac{d + K \psi}{d}\right) + K \psi \log (d + K\psi).$$
Replacing $d$, we get
$$\log L \approx K \psi \log (N \psi) + C' \psi ,$$
which yields the result. $\Box$

This is the largest quantity used for computing power sums, before
going back to the coefficients of the polynomial $\Phi[f(N\tau)]$,
whose sizes have the same order. We infer
\begin{corollary}\label{hPhi}
The height of $\Phi[f(N\tau)]$ is approximately $K
\psi(N)\log(N \psi(N))$.
\end{corollary}

%%%%% SSS
\subsubsection{The prime case}

In this section, we keep the notations $f$, $w$, etc.
\begin{proposition}\label{rootsPhi}
When $N = \ell$ is prime, the roots of $\Phi[f(\ell\tau)]$ are
$$f(\ell\tau) \text{ and } \ell^{-w} f((\tau+h)/\ell) \text{ for } 0
\leq h < \ell.$$
\end{proposition}

\medskip
\noindent
{\em Proof:} In case $N = \ell$, the system of cosets is given in
Proposition~\ref{prop:cosets}.

When $c > 0$, let $\overline{c} = 1/c \bmod \ell$, so that $R_c'
\sim_\Gamma \mat{1}{\overline{c}}{0}{\ell}$ and
$$g |_w R_c' (\tau) = \ell^{-w} f\left(\frac{\tau +
\overline{c}}{\ell}\right).$$
When $c = 0$, $R_0' = \mat{\ell}{0}{0}{1}$ is already reduced
and $g |_w R_0' (\tau) = f(\ell\tau)$.
Finally, $R_\ell' \sim_\Gamma \mat{1}{0}{0}{\ell}$, from which
$$g |_w R_\ell' (\tau) = \ell^{-w} f\left(\frac{\tau}{\ell}\right).
\quad \Box$$
In practice, though, we prefer to scale all the conjugates by
$\ell^w$, so that we have to deal with integer coefficients in the
expansions. By Corollary~\ref{hPhi}, the corresponding modular
polynomial has height close to $2 w (\ell+1) \log \ell$.

\begin{proposition}\label{Stfq}
For $t \geq 1$, the power sum $S_t(f)$ has a $q$-expansion.
\end{proposition}

\medskip
\noindent
{\em Proof:} with the notations of Proposition~\ref{propSDf}:
$$S_t(f) = \ell^{wt} f^t(\ell\tau) + S_{\ell, t}(f)$$
and both terms have a $q$-expansion. $\Box$

%%%%% SSS
\subsubsection{The case of $F_\ell$}

Remember that $\Fell{\ell}(\tau)
= E_2(\tau) - \ell E_2(\ell\tau)$ is a modular form of weight $2$ for
$\Gamma_0(\ell)$. We simplify the presentation using $F = \Fell{\ell}$
and replace $|_2$ by $|$.
We need to compute all $F|(R\tau)$ for all $R$'s from
Proposition~\ref{prop:cosets}.

\begin{proposition}\label{rootsFell}
We have $F| (R_0 \tau) = F(\tau)$. With the notations of Lemma~\ref{cbar},
for $0 < c < \ell$:
$$F| (R_c \tau) = -\frac{1}{\ell} \; F\left(\frac{\tau +
\overline{c}}{\ell}\right).$$
Also
$$F| (R_\ell\tau) = -\frac{1}{\ell} \; F\left(\frac{\tau}{\ell}\right).$$
\end{proposition}

\medskip
\noindent
{\em Proof:}
Let us suppose that $c > 0$. Start with
$$E_2\left(\mat{1}{0}{c}{1}\right)(\tau)
=(c\tau+1)^2 E_2(\tau) + \frac{6 c}{\pi i} (c \tau + 1).$$
Plugging $R_c'$ in (\ref{eqE2}), we find
$$E_2(\ell (R_c\tau)) = E_2(T_c (R_c' \tau)) =
(c (R_c'\tau)+u)^2 E_2(R_c' \tau) + \frac{6 c}{\pi i} (c (R_c' \tau)
+ u).$$
We simplify
$$c (R_c'\tau) + u = \frac{c\tau+c \overline{c} + \ell u}{\ell}
= \frac{c\tau+1}{\ell}$$
which leads to
$$E_2(\ell (R_c\tau)) = \left(\frac{c\tau+1}{\ell}\right)^2
 E_2\left(\frac{c\tau+1}{\ell}\right) + \frac{6
c}{\pi i} \frac{c\tau+1}{\ell}.$$
We deduce that
$$F(R_c\tau) = E_2(R_c\tau)-\ell E_2(\ell (R_c\tau))
= (c\tau+1)^2 \left(E_2(\tau) - \frac{1}{\ell} E_2(R_c' \tau)\right).$$
Remark that
$$F(R_c' \tau) = E_2(R_c'\tau) - \ell E_2(\tau + \overline{c}) =
E_2(R_c'\tau) - \ell E_2(\tau)$$
so that 
$$F(R_c\tau) = -\frac{1}{\ell} (c\tau+1)^2 \; F\left(\frac{\tau +
\overline{c}}{\ell}\right),$$
or
$$F| (R_c\tau) = -\frac{1}{\ell} \; F\left(\frac{\tau +
\overline{c}}{\ell}\right).$$

The last case is that of $R_\ell$, which needs
$$E_2\left(\mat{0}{-1}{1}{\ell}\right)(\tau)
=(\tau+\ell)^2 E_2(\tau) + \frac{6}{\pi i} (\tau + \ell).$$
Performing computations the way we treated $c > 0$, we find
$$F| (R_\ell \tau) = -\frac{1}{\ell} \;
F\left(\frac{\tau}{\ell}\right). \quad \Box$$

\medskip
Again, it is more convenient to scale all conjugates by $-\ell$, which
gives us $-\ell F(\tau)$ and $F((\tau+h)/\ell)$ for all $0 \leq h < \ell$.
Remark that Proposition~\ref{Stfq} applies too, {\it mutatis mutandis}.

%%%%% S
\section{Fast numerical evaluation of Eisenstein series}
\label{sct:fast}

Since one of the methods for computing modular polynomials uses
floating point evaluations, we give some algorithms to compute our
functions.

%%%%%%%%%% SS
\subsection{Jacobi $\theta$ functions}

The classical $\theta$ functions are:
$$\theta_2(\qt) = \sum_{n\in\Z} \qt^{(n+1/2)^2},
\;\theta_3(\qt) = \sum_{n\in\Z} \qt^{n^2},
\;\theta_4(\qt) = \sum_{n\in\Z} (-1)^n \qt^{n^2}.$$
Among many properties, one has
$$\theta_{3}^4(\qt) = \theta_{4}^4(\qt) + \theta_{2}^4(\qt).$$
The latter formula enables to concentrate on the evaluation of
$\theta_{3, 4}(\qt)$ as is done in \cite{Dupont11}.

Note also the following \cite[Prop. 4]{Dupont11}
\begin{proposition}\label{thinf}
$$\lim_{\mathrm{Im}(\tau) \rightarrow +\infty} \theta_{3}(\tau) = 1,
\; \lim_{\mathrm{Im}(\tau) \rightarrow +\infty} \theta_{4}(\tau) = 1,
\; \lim_{\mathrm{Im}(\tau) \rightarrow +\infty} \theta_{2}(\tau) = 0.$$
\end{proposition}

The quantities (see \cite[\S 13.20]{Erdelyi53})
$$a = \theta_{2}(\qt), \quad b = \theta_{3}(\qt), \quad c = \theta_{4}(\qt)$$
satisfy the following identities (among others)
\begin{equation}\label{E46D}
E_4 = (a^8+b^8+c^8)/2, \; E_6 = (a+b) (b+c) (c-a)/2, \; \Delta = (a b
c/2)^8.
\end{equation}
From which we deduce
$$\lim_{\mathrm{Im}(\tau) \rightarrow +\infty} E_4(\tau) = 1,
\; \lim_{\mathrm{Im}(\tau) \rightarrow +\infty} E_6(\tau) = 1.$$

%%%%%%%%%% SS
\subsection{Fast evaluation of $E_{2k}(q)$ for $k \geq 1$}

The $\theta$ functions that can be evaluated at precision $N$ in
time $O({\sf M}(N) \sqrt{N})$ with $\qt$-expansions (see
\cite{EnHaJo18}) or faster in $O({\sf M}(N)
\log N)$ using \cite{Dupont11} and also \cite{Labrande18}.
It follows that the quantities $E_{2k}$ (for $k \geq 2$)
can be evaluated at precision $N$ in $O({\sf M}(N)\log N)$
operations. As a consequence $j(q)$ can also be evaluated with the
same complexity. This is also the case for $\eta$; in practice the
lacunary properties of powers of $\eta$ can also be used (see
\cite{Serre85} and \cite[Remark 2.1.27 with the indications
therein]{CoSt17}).

\newcommand{\hypergeom}[4]{{}_2F_1\left({#1}, {#2}; {#3}; {#4}\right)}
Evaluating $E_2$ is less obvious. However, hidden in the proof of
\cite[Thm 4]{KaZa98} (thanks to \cite{KaKo03} for highlighting this),
we find
$$\frac{E_2 E_4}{E_6} = 
\frac{\hypergeom{\frac{13}{12}}{\frac{5}{12}}{1}{\frac{1728}{j}}}
     {\hypergeom{\frac{1}{12}}{\frac{5}{12}}{1}{\frac{1728}{j}}}
   = 1 + \frac{720}{j} + \cdots$$
where the Gauss hypergeometric function is defined by
$$\hypergeom{a}{b}{c}{x} = \sum_{k=0}^\infty \frac{(a)_k (b)_k}{(c)_k \,
k!} x^k, |x| < 1$$
where $(a)_k = a (a+1) \cdots (a+k-1)$. By
   \cite{Hoeven99,Hoeven01,MeSa10} and also \cite{BrZi10},
   this function can be computed at precision $N$ in $O({\sf
   M}(N)(\log N)^2)$ operations. See also \cite{Johansson19} for
   realistic computations.
Other links with hypergeometric
functions could be investigated (A.~Bostan, personal communication).

Also, note that evaluating $\Fell{\ell}$ for small prime
$\ell$ can be done using the special formulas we mentioned above. 

%%%%%%%%%% SS
\subsection{A multi-value approach}

In practice, a simpler approach yields the values $E_{2k}$ of many
$k$'s with $k \geq 1$ in time $O({\sf M}(N) \sqrt{N})$ based on
\cite{EnHaJo18}. The cost reduces to that of one series evaluation.

From \cite{BeYe02}, we take
$$(q; q)_\infty = \exp(-2i\pi\tau/24) \eta(q),$$
and for $k \geq 0$:
$$T_{2k}(q) = 1 + \sum_{n=1}^\infty (-1)^n \left\{(6n-1)^{2k} q^{n
(3n-1)/2}+ (6n+1)^{2k} q^{n (3n+1)/2}\right\}.$$
Note that $(q; q)_\infty = T_0(q)$.

\begin{theorem}[Section 6, general formulas for $T_{2k}(q)/T_0(q)$ are also
given]\label{ET}
$$\frac{T_2(q)}{T_0(q)} = E_2,
\;\frac{T_4(q)}{T_0(q)} = 3 E_2^2-2 E_4,
\;\frac{T_6(q)}{T_0(q)} = 15 E_2^3 - 30 E_2 E_4 + 16 E_6.$$
\end{theorem}

If we need to compute $E_2$, $E_4$ and $E_6$, we see that it is enough to
evaluate the series $T_{2k}$ for $k \in \{0, 1, 2, 3\}$ followed by a
handful of multiplications and divisions as given in the preceding
Theorem. Moreover, we can evaluate these series by
sharing the common powers of $q$. These powers are evaluated at a
reduced cost using \cite[Algorithm2]{EnHaJo18}. We give the modified
procedure as algorithm~\ref{algo1}.
In Step 3.2.3, we have added the contribution $(6n \pm 1)^{2i}$ to each
\verb+T[i]+. We assume that the cost of multiplying by these small
quantities is negligible. Were it not the case, we could use
incremental computations of the polynomials $(6n \pm 1)^{2i}$. The
cost of this algorithm reduces to that of one of the series, gaining a
factor $kmax$.

\LinesNotNumbered
\begin{algorithm}[hbt]
\caption{Combined evaluation of $T_{2k}(q)$. \label{algo1}}
\SetKwProg{Fn}{Function}{}{}
\Fn{EvaluateManyT(q, N, kmax)}{
\Input{$q$, $N$, $kmax$}
\Output{$(T_{2k}(q))$ for $0 \leq k \leq kmax$}

1. \For{$k:=0$ \KwTo $kmax$}{
    $T[k] \leftarrow 0$\;
}

2. $s \leftarrow 1$; $A \leftarrow \{1\}$; $Q[1] \leftarrow \{q\}$; $c
\leftarrow 0$\;

3. \For{$n:=1$ while $n (3n+1)/2 \leq N$}{
$s \leftarrow -s$\; \tcp*[h]{$s = (-1)^n$}

3.1 $c \leftarrow c + 2 n-1$\;

3.2 \For{$r:=1$ \KwTo $2$}{
 3.2.1 \If{$r = 2$}{
    $c \leftarrow c + n$\; \tcp*[h]{$c = n (3n+1)/2$}
 }
 3.2.2 $q' \leftarrow $ FindPowerInTable($A$, $Q$, $c$)\;
 3.2.3 $C \leftarrow (6n+(-1)^r)^2$\;
 3.2.4 \For{$k:=0$ \KwTo kmax}{
     $T[k] \leftarrow T[k] + s q'$\;
     \If{$k < kmax$}{
         $q' \leftarrow C q'$\;
     }
 }}}
4. \For{$k:=0$ \KwTo kmax}{
   $T[k] \leftarrow T[k] + 1$\;
   }
5. \Return{$T$}.
}
\end{algorithm}
Algorithm~\ref{algo1} uses the primitive in Algorithm~\ref{algo0}. The
reason of Step 4 is that $T[k]$ will be close to $1$ when $q$ is
small, so that we may not want to add 1 right at the beginning and
perhaps not in this function.

\LinesNotNumbered
\begin{algorithm}[hbt]
\caption{Finding $c$ as a combination of known values. \label{algo0}}
\SetKwProg{Fn}{Function}{}{}
\Fn{FindPowerInTable($A$, $Q$, $c$)}{
\Input{$A = \{a_1, \ldots, a_z\}$, $Q$ such that for all $i$, $Q[a_i]
= q^{a_i}$, $c$}
\Output{$q^c$; $A$ and $Q$ are updated}

 \If{$c = 1$}{
    $q' \leftarrow Q[1]$\;
 }
 \ElseIf{$c = 2a$ with $a \in A$}{
    $q' \leftarrow Q[a]^2$\;
 }
 \ElseIf{$c = a+b$ with $a, b \in A$}{
    $q' \leftarrow Q[a] \cdot Q[b]$;
 }
 \ElseIf{$c = 2 a+b$ with $a, b \in A$}{
    $q' \leftarrow Q[a]^2 \cdot Q[b]$\;
 }
   $A \leftarrow A \cup \{c\}$\;
   $Q[c] \leftarrow q'$\;
   \Return{$q'$}.
}
\end{algorithm}

%%%%%%%%%% SS
\subsection{The case of imaginary arguments}

In practice, it is easier to consider $\tau = \rho i$ for real $\rho \geq
1$. In that case, $1 > q_0 = \exp(-2\pi) = 0.001867\ldots \geq q =
\exp(-2\pi \rho) > 0$. The functions $E_2$ and $E_6$ are increasing
from $E_{2k}(q_0)$ to
$1$ (note that $E_6(q_0) = 0$ and $E_2(q_0) = 3/\pi$ from~\cite{ElSe10});
$E_4$ is decreasing from $E_4(q_0)$ to
$1$. This is important to note for the computations not to
explode. Remember also that $j(i) = 1728$.

We turn to the precision needed for evaluating the functions $T_{2k}$.
Let $N$ denote an integer and $T_{2k, N}$ the truncated sum up to $n =
N-1$. Since the series is alternating, we can bound the error using
\begin{eqnarray*}
|T_{2k}(q)-T_{2k, N}(q)| & \leq & \{(6N-1)^{2k} \, q^{N(3 N-1)/2} +
(6N+1)^{2k} \, q^{N(3 N+1)/2}\} \\
&\leq & ((6N-1)^{2k}+(6N+1)^{2k}) \, q^{N(3 N-1)/2}.
\end{eqnarray*}

\noindent
Since $0 < q < 1$, this gives us a very fast quadratic convergent series.

%%%%% S
\section{The polynomials of Fricke and Charlap/Coley/Robbins}
\label{sct:Fricke}

%%%%%%%%%% SS
\subsection{The work of Elkies}

An isogeny is associated with its kernel, or its polynomial
description (called {\em kernel polynomial}). Given some finite
subgroup $F$ of $\EE$, one can build an
isogenous curve $\Es$ and the corresponding isogeny, using V\'elu's formulas.
In the context of point counting, we discover
a curve $\Es$ that is $\ell$-isogenous to $\EE$ via its $j$-invariant
as a root of the traditional
modular polynomial, and we need to find the coefficients of $\Es$,
together with the isogeny. The idea of Elkies is to consider the same
problems on the Tate curves associated to the elliptic curves $\EE$
and $\Es$.

To be brief, $\EE$ has an equation in some parameter $q$, and the
isogenous $\Es$ is associated to parameter $q^\ell$, where $\ell$ is
the degree of the isogeny, which in our case is associated with a
finite subgroup $F$ of cardinality $\ell$. To be more precise, we
consider $\EE$ has having equation $y^2 = x^3 + A x + B$ with
\begin{equation}\label{ABE4E6}
A = -3 E_4(q), B = -2 E_6(q).
\end{equation}
With a compatible scaling, we get the equation for $\Es: y^2 = x^3 +
\As x + \Bs$ with
\begin{equation}\label{AsBsE4E6}
\As = -3 \ell^4 E_4(q^\ell), \; \Bs = -2 \ell^6 E_6(q^\ell).
\end{equation}
More importantly, writing $\kappa_r$ for the power sums of the roots of
the kernel polynomial, we have
\begin{equation}\label{eqsig}
\kappa_1 = \frac{\ell}{2} (\ell E_2(q^\ell)-E_2(q)) =
-\frac{\ell}{2} \Fell{\ell}(q).
\end{equation}
Beyond this, Elkies proved~\cite[formulas (66) to (69)]{Elkies98}
\begin{proposition}\label{Elkies:formulas}
$$A-\As = 5 (6 \kappa_2 + 2 A \kappa_0),$$
$$B-\Bs = 7 (10 \kappa_3 + 6 A \kappa_1+4 B \kappa_0),$$
together with an induction relation satisfied by other $\kappa_k$ for
$k > 3$.
\end{proposition}
This can rephrased as $(\kappa_1, \As, \Bs)$ is enough to describe an
isogeny. Also $\As$ and $\Bs$ belong to $\Q[\kappa_1, A, B]$
since $\kappa_2$ and $\kappa_3$ do. The minimal polynomial of
$2 \kappa_1$ is the modular polynomial associated to $\Fell{\ell}$,
and we can express $\As$ and $\Bs$
as elements in the field $\Q[\kappa_1, A, B]$, which we use
below. Rephrased another times, $E_4(q^\ell)$ and $E_6(q^\ell)$ are
modular forms we need to express as expressions in known modular
forms. See~\cite{Elkies98} for more details on this subject.

Given these quantities, there are several algorithms to get the
isogeny. We refer to \cite{BoMoSaSc08} for this.

%%%%%%%%%% SS
\subsection{The Fricke polynomials}

%%%%%%%%%%%%%%% SSS
\subsubsection{Reinterpreting Elkies's results}

One way of looking at the work of Elkies (taken from
\cite{ChCoRo91}, but which originated in~\cite{Fricke22}) is to
realize that we try to decompose the $\ell$-th division polynomial
$f_\ell$ (say $\ell$ is odd) of degree $(\ell^2-1)/2$ over a subfield
of degree $\ell+1$. In Fricke's term, we compute a degree $\ell+1$
resolvent for the equation $f_\ell(X) = 0$.

\begin{center}
\begin{tikzpicture}[scale=0.75]
\node (a) at (2, 2) {$\Q(A, B)[X]/(f_\ell(X, A, B))$};
\node (b) at (2, 0) {$\Q(A, B)[X]/(U_\ell(X, A, B))$};
\node (c) at (2, -2) {$\Q(A, B)$};
\draw[-] (a) -- (b) node[midway,right] {$(\ell-1)/2$};
\draw[-] (b) -- (c) node[midway,right] {$\ell+1$};
\end{tikzpicture}
\end{center}

%%%%%%%%%% SS
\subsubsection{Theory}

We start from an elliptic curve $\EE: y^2 = x^3+A x+B$ and we fix some
odd prime $\ell$, putting $d = (\ell-1)/2$. Our aim is to find the
equation of an $\ell$-isogenous curve $\Es: y^2 = x^3 + \As x +
\Bs$. The results for $U_\ell$ are due to Fricke, and
Charlap/Coley/Robbins for $V_\ell$ and $W_\ell$.

\begin{theorem}
There exist three polynomials $U_\ell$, $V_\ell$,
$W_\ell$ in $\Z[X, Y, Z, 1/\ell]$ of degree $\ell+1$ in $X$
such that $U_\ell(2 \kappa_1, A, B)=0$, respectively $V_\ell(\As, A, B)
= 0$, $W_\ell(\Bs, A, B) = 0$.
\end{theorem}

Let us turn our attention to the properties of these polynomials.
Note that $U_\ell$ (resp. $V_\ell$
and $W_\ell$) is the minimal polynomial of a weight 2 form (resp. 4
and 6). We note $\varpi$ the corresponding weight. 
\begin{theorem}
When $\ell > 3$, the polynomials $U_\ell$, $V_\ell$, $W_\ell$ live in
$\Z[X, Y, Z]$.
\end{theorem}

As a consequence of Theorem~\ref{mp1}, we have
\begin{proposition}
The polynomials $U_\ell$, $V_\ell$ and $W_\ell$ are homogeneous with
weight $\varpi (\ell+1)$.
\end{proposition}

\begin{proposition}\label{prop5.5}
Put $z^\ell = q = \exp(2 i \pi \tau)$ and $\zeta_\ell$ a root of
unity. Then

1) The roots of $U_\ell(X, A(q), B(q))$ are 
$\Fell{\ell}(z \zeta_\ell^h)$ for $0 \leq h < \ell$, and
$-\ell \Fell{\ell}(q)$.

2) The roots of $V_\ell(X, E_4(q), E_6(q))$ (resp. $W_\ell$) are
$E_4(z \zeta_\ell^h)$ (resp.
$E_6(z \zeta_\ell^h)$) for $0
\leq h < \ell$, and $\ell^4 E_4(q^\ell)$ (resp. $\ell^6
E_6(q^{\ell}$)).
\end{proposition}

Part 1 is proven in Proposition~\ref{rootsFell}; part 2 is
done in Proposition~\ref{rootsPhi}. We can also use
Proposition~\ref{sizeME} to get
\begin{proposition}\label{sizeUVW}
The height of $U_\ell$ (resp. $V_\ell$, $W_\ell$) is approximately
$\varpi (\ell+1) \log\ell$.
\end{proposition}

%%%%%%%%%%% SSS
\subsubsection{Representing $\As$ and $\Bs$ as rational fractions}
\label{sct:rat}

Once we have computed $U_\ell$, we can either compute $V_\ell$ and
$W_\ell$ or use another representation.
From \cite[Theorem 3.9]{NoYaYo20}, there exist polynomials
$\mathcal{A}_\ell$ and $\mathcal{B}_\ell$ of degree less than $\ell+1$
such that
\begin{equation}\label{eqAs0}
\As = \frac{\mathcal{A}_\ell(X, A, B)}{U_\ell'(X)}, \; \Bs =
\frac{\mathcal{B}_\ell(X, A, B)}{U_\ell'(X)}
\end{equation}
(Only here: $U_\ell'(X) = \frac{\partial U_\ell}{\partial X}$.)
Moreover, $\mathcal{A}_\ell$ and $\mathcal{B}_\ell$ are polynomials
with integer coefficients and of respected generalized weight $2 \ell+4$
and $2\ell+6$. These formulas are of independent interest and may
prove useful in other contexts.
The authors of the reference use Groebner basis
computations to find the two numerators. We propose another route
later on.

%%%%%%%%%% SSS
\subsubsection{Computing isogenous curves over finite fields}
\label{sct:summary}

When using $(U_\ell$, $V_\ell$, $W_\ell)$, we need to find the roots of
three polynomials of degree $\ell+1$ instead of a single one in the
traditional case. In general, if
$U_\ell$ has rational roots (it should be 1, $2$ or $\ell+1$), then
this is the case for each of $V_\ell$, $W_\ell$. For each triplet of
solutions $(2 \kappa_1, z_1, z_2)$ we need to test whether this leads to
an isogeny or not. See techniques for this task in \cite{BoMoSaSc08}.
Using the rational fractions for $\As$ and $\Bs$ is faster, just
needing evaluations of rational fractions. A method based on using
$U_\ell$ only is described in Section~\ref{sct:FrickeA}.

%%%%% S
\section{Computing Fricke polynomials}
\label{sct:computations}

The authors of \cite{ChCoRo91} give two methods of computation using
manipulations of $q$-expansions of series over
$\Q$. Following Atkin~\cite{Atkin88b}, these can be replaced by
computations modulo small primes (preferably with convenient FFT
multiplication) followed by recovery using the Chinese remaindering
theorem using the bounds in Proposition~\ref{sizeUVW}. To this, we add
two evaluation-interpolation algorithms already used for classical modular
polynomials: the first is based on floating point calculations, the
second on isogeny volcanoes.

The polynomials $U_\ell$, $V_\ell$ and $W_\ell$ are modular
polynomials of modular forms, so that they obey Theorem~\ref{mp1}. Also
$$\mathcal{A}_\ell(X, Y, Z) = \sum_{r=0}^\ell X^r \sum_{2i_2+3i_3 = \ell+2-r}
a_{r, i_2, i_3} Y^{i_2} Z^{i_3},$$
$$\mathcal{B}_\ell(X, Y, Z) = \sum_{r=0}^\ell X^r \sum_{2i_2+3i_3 = \ell+3-r}
b_{r, i_2, i_3} Y^{i_2} Z^{i_3}.$$
All the methods to be described can be applied to $U_\ell$, $V_\ell$,
$W_\ell$, and also to $\mathcal{A}_\ell$,
$\mathcal{B}_\ell$. To simplify the presentation, we assume from now
on (unless indicated) that $\ell > 3$ and concentrate on $U_\ell$,
indicating what has to be changed for the other polynomials; in
particular we assume we are looking for the modular polynomial of a
form of weight $w$. 

We rewrite
$$U_\ell(X) = X^{\ell+1} + C_1(E_4, E_6, \Delta) X^{\ell} + \cdots +
C_{\ell+1}(E_4, E_6, \Delta).$$
By Theorem~\ref{mp1}, $C_t$ is a modular form for $\Gamma$ of weight
$w t = 2t$; this implies $C_1 = 0$. By (\ref{Eeq}), $C_2 = c_2 E_4$,
$C_3 = c_3 E_6$, $C_4 = c_4 E_4^2$, $C_5 = c_5 E_4 E_6$, $C_6 = c_{6,
1} E_4^3 + c_{6,0} \Delta$. For instance, the following methods will
give us 
$$U_5(X) = X^6-60 E_4 X^4-320 E_6 X^3-720 E_4^2 X^2-768 E_4 E_6 X-320
E_4^3+552960 \Delta.$$

%%%%%%%%%%%%%%%	SSS
\subsection{Using $q$-expansions}
\label{sct:qexp}

Note that
$$2 \kappa_1(q) = -\ell\Fell{\ell}(q) = \ell (\ell-1) + 24 \ell
\sum_{n=1}^\infty \delta_1'(n) q^n$$
where $\delta_1'(n)$ is the sum of the divisors of $n$ prime to $\ell$.

Using Proposition~\ref{prop5.5}, we denote by
$\sigma_t(q)$ the corresponding power sums of roots of $U_\ell$:
$$\sigma_t(q) = (-\ell \Fell{\ell}(q))^t + \sum_{h=0}^{\ell-1}
\Fell{\ell}(z \zeta_\ell^h)^t$$
for $z^\ell = q$ and $1 \leq t \leq \ell+1$. We compute them and
recover the coefficients of $U_\ell$ using Newton's formulas as
explained in Section~\ref{sct:mpmf}.

The power sums $\sigma_t(q)$ are modular forms of weight
$w t$ and can be represented as polynomials in $E_4$, $E_6$, $\Delta$
$$\sigma_t(q) = \sum_j c_{t, j} P_{w t, j}$$
using Proposition~\ref{prop1} and Algorithm~\ref{algoAllP}.
This leads to a triangular linear system $\mathcal{S}_t$ in the $c_{t,
j}$'s. The system has $\approx (w t/6)$ rows and can be solved with
$O((w t)^{2})$ operations over $\Z$, for a total of $O(\sum_t (w t)^{2}) =
O(\ell^{3})$. Once solved for all $t$'s, we use
Newton's identities to recover the coefficients of $U_\ell$.

To start the process, one needs to evaluate the series $\sigma_t(q)$
using intermediate expressions in $z = q^{1/\ell}$
having roots of unity $\zeta_\ell$ temporarily appearing and
vanishing. See \cite[\S 2.2]{Enge09b} for more details and complexity analysis.
In particular, if we denote by $\sfM_q(d)$ the number of arithmetic
operations in $\Z$ required to multiply two dense $q$-expansions with
$d$ terms, then the total complexity of the series computations is
$O(\ell \sfM_q(\ell d))$, which is $O(\ell \sfM_q(\ell^2))$ in our
case. If $H$ is a bound on the height of the polynomial, then the bit
complexity is $O(\ell^3 (\log \ell) \sfM(H))$. Assuming $H \in
O(\ell\log\ell)$ by Proposition~\ref{sizeUVW}, this is $O(\ell^4
\log^{3+\epsilon} \ell)$.

\medskip
\noindent
{\bf Example.} Consider the case $\ell = 5$. The systems
$\mathcal{S}_t$ to be solved come from the equations:
\begin{eqnarray*}
\sigma_2(q) &=& c_{2, 0} E_4, \\
\sigma_3(q) &=& c_{3, 0} E_6, \\
\sigma_4(q) &=& c_{4, 0} E_4^2, \\
\sigma_5(q) &=& c_{5, 0} E_6 E_4, \\
\sigma_6(q) &=& c_{6, 0} E_4^3 + c_{6, 1} \Delta.
\end{eqnarray*}
We compute
$$\sigma_6(q) = 1000320+186071040 q+ \cdots$$
and we remember that $E_4(q) = 1+240 q + \cdots$, $\Delta(q) = q + \cdots$ so
that the system $\mathcal{S}_6$ is
$$\left\{\begin{array}{rcl}
1000320 &=& c_{6, 0} \\
186071040 &=& 720 \, c_{6, 0} + c_{6, 1}\\
\end{array}\right.$$
which is triangular indeed and therefore easy to solve. Its solutions
are integers.

We can also work modulo small primes and use the Chinese Remainder
Theorem to recover the polynomials. 

%%%%%%%%%% SS
\subsection{Floating point methods}

We adapt the methods proposed for ordinary modular equations to our
polynomials $(U_\ell, V_\ell, W_\ell)$. We note $H$ for the
logarithmic height of the polynomials, that we have estimated to $\varpi
(\ell+1)\log\ell$ in Proposition~\ref{sizeUVW}. All the methods are
heuristic.

%%%%%%%%%%%%%%% SSS
\subsubsection{Solving a linear system}

We start from
$U_\ell(2 \kappa_1(q), E_4(q), E_6(q), \Delta(q)) = 0$
and we compute floating point values to get a linear system in the
coefficients that should come out as integers for $\ell > 3$. We
evaluate $\kappa_1(q)$, $E_4(q)$ and $E_6(q)$ (and therefore
$\Delta(q)$ at high precision for chosen values of (imaginary) $\tau$
in $q = \exp(2i\pi\tau)$. This would involve $O(\ell^{2\omega})$
floating point operations, and we can do better in the following section.

%%%%%%%%%%%%%%% SSS
\subsubsection{Using power sums}

The case of the traditional modular polynomial is treated in
\cite{Enge09b}. We can use the same approach for our polynomials.
First of all, we need to compute
$$-{\ell} \Fell{\ell}(q), \quad \{\Fell{\ell}(z \zeta_\ell^h), 0 \leq
h < \ell\},$$
where $z^{\ell} = q$ and $\zeta_\ell$ is a primitive $\ell$-th root of
unity. By definition
$$\Fell{\ell}(q) = E_2(q)-\ell E_2(q^\ell)$$
and
$$\Fell{\ell}(z \zeta_\ell^h) = E_2(z \zeta_\ell^h) - \ell
E_2(q),$$
and the last term is a constant w.r.t. $h$. We first evaluate
$E_2(q)$, $E_2(q^\ell)$ and then the other roots, by sharing the
computations: All terms we need are of the form $(z \zeta_\ell^h)^e = z^e
\zeta_\ell^{(h e) \bmod \ell}$. When $\ell \mid e$, the computation is
a little faster. We give the corresponding code as
Algorithm~\ref{algo2}. We also precompute $\xi_h = \zeta_\ell^h$. The
complete code is in Algorithm~\ref{algo3}.
Multiple evaluations of $\eta(k \tau)$ can be shared as explained
in~\cite{EnHaJo18}.

\LinesNotNumbered
\begin{algorithm}[hbt]
\caption{Combined evaluation of $T_{2k}(z \zeta_\ell^h)$. \label{algo2}}
\SetKwProg{Fn}{Function}{}{}
\Fn{EvaluateConjugateValues($\ell$, $z$, $N$, $(\xi)$, $kmax$)}{
\Input{$\ell$, $z$, $(\xi)$, $N$}
\Output{$(T_{2k}(z \zeta_\ell^h))$ for $0 \leq k \leq kmax$, $0 \leq h
< \ell$}

1. \For{$k:=0$ \KwTo $kmax$}{
     \For{$h:=0$ \KwTo $\ell-1$}{
         $T[k, h] \leftarrow 0$\;
     }
}

2. $s \leftarrow 1$; $A \leftarrow \{1\}$; $Z[1] \leftarrow \{z\}$; $c
\leftarrow 0$\;

3. \For{$n:=1$ while $n (3n+1)/2 \leq N$}{
$s \leftarrow -s$\; \tcp*[h]{$s = (-1)^n$}

3.1 $c \leftarrow c + 2 n-1$\;

3.2 \For{$r:=1$ \KwTo $2$}{
 3.2.1 \If{$r = 2$}{
    $c \leftarrow c + n$\; \tcp*[h]{$c = n (3n+1)/2$}
 }
 3.2.2 $z' \leftarrow s \cdot$FindPowerInTable($A$, $Z$, $c$)\;
 3.2.3 $C \leftarrow (6 n+(-1)^r)^2$\;
 3.2.4 \For{$k:=0$ \KwTo kmax}{
     $T[k, 0] \leftarrow T[k, 0] + z'$\;
     \For{$h:=1$ \KwTo $\ell-1$}{
         $T[k, h] \leftarrow T[k, h] + \xi[(h c) \bmod \ell] \, z'$\;
     }
     \If{$k < kmax$}{
         $z' \leftarrow C \cdot z'$\;
     }
}}}
4. \For{$k:=0$ \KwTo $kmax$}{
     \For{$h:=0$ \KwTo $\ell-1$}{
         $T[k, h] \leftarrow T[k, h] + 1$\;
     }
   }
5. \Return{$T$}.
}
\end{algorithm}

\LinesNotNumbered
\begin{algorithm}[hbt]
\caption{Computing Fricke polynomial using floating point
 numbers. \label{algo3}}
\SetKwProg{Fn}{Function}{}{}
\Fn{ComputeUVW($\ell$, $w$)}{
\Input{$w \in \{2, 4, 6\}$ corresponding to $U_\ell$, $V_\ell$ or
$W_\ell$ respectively, $\ell$ an odd prime}
\Output{the corresponding Fricke polynomial}

0.0 $H \leftarrow w (\ell+1) \log \ell$; all computations are carried
out at precision $H$\;

0.1 compute $\zeta_\ell \leftarrow \exp(2 i \pi/\ell)$\;

0.2 \For{$h:=0$ \KwTo $\ell-1$}{
    $\xi_h \leftarrow \zeta_\ell^h$\;
   }

0.3 Compute all systems $\mathcal{S}_t$ for $2 \leq t \leq \ell+1$\;

0.4 \For{$t:=2$ \KwTo $\ell+1$}{
    $\mathcal{L}_t \leftarrow \emptyset$\;
    }

0.5 $\rho \leftarrow 1$\;

1. \While{there is a system $\mathcal{L}_t$ that is not solved}{

1.0 $\rho \leftarrow \rho + 0.1$\;

1.1 $z \leftarrow \exp(-2 \pi \rho/\ell)$; $q_\rho \leftarrow z^\ell$\;

1.2 $T \leftarrow $ EvaluateConjugateValues($\ell$, $z$, $N$, $(\xi_h)$,
$s$)\;

1.3. use Theorem~\ref{ET} to evaluate $E_{2k}$ for all $q_i$'s from $T$,
yielding $(\kappa(q_i))$ for $i = 0, \ldots, \ell+1$; also deduce
$E_{4, \rho} = E_4(q_\rho)$, $E_{6, \rho} = E_6(q_\rho)$,
$\Delta_{\rho} = (E_{4, \rho}^3-E_{6, \rho}^2)/1728$\;

1.4 \For{$r:=t$ \KwTo $\ell+1$}{
  \If{$\mathcal{L}_t$ is not solved}{

1.4.1 Instantiate $\mathcal{S}_t$ with $\sum_i \kappa_1(q_i)^t$,
  $E_{4, \rho}$, $E_{6, \rho}$, $\Delta_\rho$; add it to $\mathcal{L}_t$\;

1.4.2 \If{$\mathcal{L}_t$ has as many equations as unknowns}{
   solve $\mathcal{L}_t$ and store the values; declare $\mathcal{L}_t$ solved\;
  }

  }
}
}
2. Round the coefficients and use Newton's formulas.
}
\end{algorithm}
The system $\mathcal{S}_t$ has size $O(t^2)$ and we
need $O(t^\omega)$ operations to solve it, for a total of
$O(\ell^{\omega+1})$. Like in the series case, anticipate integer
coefficients, which makes recognition of the coefficients easier.

\smallskip
\noindent
{\bf Example.} Take again $\ell=5$, for which the $\mathcal{S}_t$ were
already given. Let us concentrate on the case of
$\sigma_6(q) = u_{6, 0} E_4^3 + u_{6, 1} \Delta$; we start with
$\mathcal{L}_6 = \emptyset$. Using $\rho = 1.1$ leads to
$$\mathcal{L}_6 = \{
1.912407642 u_{6,0}+0.0009726854527956 u_{6,1}
= 1393450.57337539139\}$$
and the following iteration with $\rho=1.2$ adds
$$\{
1.435895343 u_{6,0}+0.0005247501300701 u_{6,1} 
= 1156054.63606077432\}$$
and the solution of $\mathcal{L}_6$ (rounded to integers) is
$$u_{6, 1} = -534159360, u_{6, 0} = 1000320.$$

%%%%%%%%%% SS
\subsection{Isogeny volcanoes}

The method in \cite{ChLa05} shares many common points
with the method to be described next but
with a worse complexity. It uses supersingular curves whose
complete explicit $\ell$-torsion is required.
The work of \cite{BrLaSu12} is a building block in \cite{Sutherland13}
where direct evaluation of $\Phi_\ell(X, j(E)) \bmod q$ is made
possible using an explicit version of the Chinese remainder theorem
modulo small primes. Our version is an adaptation to the computation
of the Fricke polynomials.

%%%%%%%%%%%%%%% SSS
\subsubsection{Quick presentation}

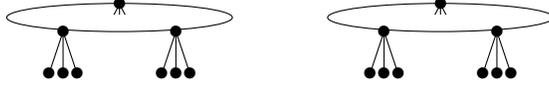
\begin{figure}[hbt]
\begin{center}
\begin{tikzpicture}[scale=0.75]
\draw (0, 0) ellipse [x radius=2, y radius=0.25];
%%%%% first tree
\node at (-1, -0.25) {$\bullet$};
 \draw[-] (-1, -0.25) -- (-1.25, -1); \node at (-1.25, -1) {$\bullet$};
 \draw[-] (-1, -0.25) -- (-1, -1); \node at (-1, -1) {$\bullet$};
 \draw[-] (-1, -0.25) -- (-0.75, -1); \node at (-0.75, -1) {$\bullet$};
%%%%% second tree
\node at (1, -0.25) {$\bullet$};
 \draw[-] (1, -0.25) -- (1.25, -1); \node at (1.25, -1) {$\bullet$};
 \draw[-] (1, -0.25) -- (1, -1); \node at (1, -1) {$\bullet$};
 \draw[-] (1, -0.25) -- (0.75, -1); \node at (0.75, -1) {$\bullet$};
%%%%% third tree beyond
\node at (0, 0.23) {$\bullet$};
 \draw[-] (0, 0.23) -- (-0.1, 0.05);
 \draw[-] (0, 0.23) -- (0, 0.05);
 \draw[-] (0, 0.23) -- (0.1, 0.05);
\end{tikzpicture}
\hspace*{1cm} \begin{tikzpicture}[scale=0.75]
\draw (0, 0) ellipse [x radius=2, y radius=0.25];
%%%%% first tree
\node at (-1, -0.25) {$\bullet$};
 \draw[-] (-1, -0.25) -- (-1.25, -1); \node at (-1.25, -1) {$\bullet$};
 \draw[-] (-1, -0.25) -- (-1, -1); \node at (-1, -1) {$\bullet$};
 \draw[-] (-1, -0.25) -- (-0.75, -1); \node at (-0.75, -1) {$\bullet$};
%%%%% second tree
\node at (1, -0.25) {$\bullet$};
 \draw[-] (1, -0.25) -- (1.25, -1); \node at (1.25, -1) {$\bullet$};
 \draw[-] (1, -0.25) -- (1, -1); \node at (1, -1) {$\bullet$};
 \draw[-] (1, -0.25) -- (0.75, -1); \node at (0.75, -1) {$\bullet$};
%%%%% third tree beyond
\node at (0, 0.23) {$\bullet$};
 \draw[-] (0, 0.23) -- (-0.1, 0.05);
 \draw[-] (0, 0.23) -- (0, 0.05);
 \draw[-] (0, 0.23) -- (0.1, 0.05);
\end{tikzpicture}
\end{center}
\caption{A typical set of volcanoes. \label{volcanoes}}
\end{figure}

In a nutshell, the algorithm in~\cite{BrLaSu12} performs computations
modulo special primes $p$ satisfying arithmetical conditions: $p \equiv
1\bmod \ell$ and $4 p = t^2 - \ell^2 v^2 D$ in integers $t$ and $v$, $D$
not a multiple of $\ell$; $D < 0$ is the (fundamental) discriminant of some
auxiliary quadratic field. With these conditions, the so-called class
polynomial $H_D(z)$ splits completely modulo $p$ and its roots are
$j$-invariants of elliptic curves with complex multiplication by
the maximal order $\mathcal{O}_D$. The isogeny volcanoes that we can
build have only one level (see Figure~\ref{volcanoes}) and the
corresponding $j$-invariants are the roots of $H_{\ell^2 D}(X)$.
Basically, the algorithm interpolates data
using the isogenies attached to the volcanoes. 
We refer the reader to the original article for more properties
related to elliptic curves. For our purpose, we just need to know that
we have isogeny data available and that they can help us computing the
polynomial $U_\ell(X, E_4, E_6, \Delta) \bmod p$ from these data. We
refer to the article for the complexity under GRH, namely $O(\ell^3
(\log \ell)^3 \log\log\ell)$ using $O(\ell^2 \log (\ell p))$ space for
suitably chosen $p$.

We adapt a slight modification of the simplified version Algorithm
2.1 of \cite{BrLaSu12} to our needs to give Algorithm~\ref{core}. All
we describe is also valid in the full version in \cite{BrLaSu12}.

\LinesNotNumbered
\begin{algorithm}[hbt]
\caption{The core algorithm. \label{core}}
\SetKwProg{Fn}{Function}{}{}
\Fn{PartialVolcano($\ell$, $D$, $H_D(X)$, $p$)}{
\Input{$\ell$ an odd prime, $D$ the discriminant of an imaginary
quadratic order $\mathcal{O}$ with class number
$h(D) \geq \ell+2$; $H_{D}$ the class
polynomial associated to the order $\mathcal{O}$; $p$ prime with
$p\equiv 1\bmod\ell$ and $4 p = t^2 - \ell^2 v^2 D$, $v\not\equiv 0
\bmod \ell$}
\Output{A collection $(\EE_i, \{Q_{ik}, \EE_{ik}'\}_{1 \leq k \leq
\ell+1})_{1\leq i\leq h}$ where 
$\EE_{ik}' = \EE_i / \langle Q_{ik}\rangle$}
1. Build the list $\mathcal{J}_D$ containing the roots of $H_{D}(z)$ 
modulo $p$\;
2. \For{$j_i \in \mathcal{J}_D$}{

2.1 find a curve $\EE_i: y^2 = x^3 + A_i x + B_i$ having invariant
$j_i$ and cardinality $m = p+1-t$;

2.2 find all the neighbors $\mathcal{N}(\EE_i)$ in the volcano of $\EE_i$: $2$
horizontal isogenies and $\ell-1$ on the floor. Let $\mathcal{K}$ be a
set to contain invariants and initialized to $\emptyset$.

\While{we do not have all isogenies of both kinds}{

2.2.1 Select a random point $Q_{ik}$ of order $\ell$
on $\EE_i/\GFq{p} = [A_i, B_i]$. 

2.2.2 Compute the rational isogeny $\EE_i
\rightarrow \EE_{ik}' = \EE_i / \langle Q_{ik}\rangle$ using V\'elu's
formulas. The result is a pair $(Q_{ik}, \EE_{ik}')$.

2.2.3 \If{$j(\EE_{ik}') \not\in \mathcal{K}$}{

$\mathcal{K} \leftarrow \mathcal{K} \cup \{j(\EE_{ik}')\}$\;

if $j(\EE_{ik}')$ is a root of $H_{D}$,
then $\EE_{ik}'$ is on the crater and is one of the two neighbours. If
it does not belong to the crater, it belongs to the floor. Store the
pair $(Q_{ik}, \EE_{ik}')$.
}
}
Store $\EE_i, \{Q_{ik}, \EE_{ik}'\}_{1 \leq k \leq \ell+1}$\;
}
3. \Return{$(\EE_i, \{Q_{ik}, \EE_{ik}'\}_{1 \leq k \leq
\ell+1})_{1\leq i\leq h}$}.
}
\end{algorithm}

%%%%%%%%%%%%%%% SSS
\subsubsection{The algorithm for $U_\ell$}

We denote by $\powersums_r$ the power sums of $U_\ell$. Such a
$\powersums_r$ is a modular form of weight $w r$ for $\Gamma$.
As explained in Section~\ref{sss:express}, we may write these power
sums as
$$\powersums_r(E_4, E_6, \Delta) = \sum_{j_r=0}^{m_r}
c_{r, j_r} P_{w r, j}.$$
The values of $P_{w r, j}$ and the sums will be reconstructed from values
$\kappa_{1, i}$ associated to curves $\EE_i: y^2 = x^3 +
A_i x+B_i$.

\LinesNotNumbered
\begin{algorithm}[hbt]
\caption{Computing $U_\ell(X, Y, Z) \bmod p$. \label{algo-U}}
\SetKwProg{Fn}{Function}{}{}
\Fn{\ComputeUMod($\ell$, $D$, $H_{D}(z)$, $p$):}{
\Input{$\ell$ an odd prime, $D$ the discriminant of an imaginary
quadratic order $\mathcal{O}$ of discriminant $D$ with class number
$h(D) \geq \ell+2$; $H_{D}$ is the class
polynomial associated to order $\mathcal{O}$; $p$ prime with
$p\equiv 1\bmod\ell$ and $4 p = t^2 - \ell^2 v^2 D$, $v\not\equiv 0
\bmod \ell$}
\Output{$U_\ell(X, Y, Z) \bmod p$}

1. $(\EE_i, \{Q_{ik}, \EE_{ik}'\}_{1 \leq k \leq
\ell+1})_{1\leq i\leq h} \leftarrow $PartialVolcano($\ell$, $D$,
$H_D(X)$, $p$)\;

2. Evaluate the quantities $P_{w t, j}$ in the $\{\EE_{ik}'\}$ using
Algorithm~\ref{algoAllP}\;

3. \For{$i \leftarrow 1$ \KwTo $\ell+1$}{

 \For{$t \leftarrow 1$ \KwTo $\ell+1$}{
 $\sigma_{t, i} \leftarrow (1/2) \sum_{k=1}^{\ell-1} x(Q_{ik})^t$\;
}
}

4. \For{$t \leftarrow 1$ \KwTo $\ell+1$}{
solve the linear system
$$\sigma_{t, i} = \sum_{j=0}^{m_t} c_{t, j} P_{t, j}.$$\;
}

5. \Return{$U_\ell$ recovered from $c_{t, j_t}$'s using Newton's
formulas.}
}
\end{algorithm}
Note that for each $i$, the kernel polynomial of the isogeny from
$\EE_i$ to $\EE_{ik}'$ starts $X^{(\ell-1)/2} - \kappa_{1, i}
X^{(\ell-3)/2} + \cdots$, so that
$$U_\ell(X, E_{4, i}, E_{6, i}, \Delta_i) = \prod_{i=1}^{\ell+1} (X -
\kappa_{1, i}).$$
Given these roots, it is easy to compute the power sums, see
Algorithm~\ref{algo-U}. Trading multiplications for additions, Step 2 costs
$O(\ell^2)$ operations over $\GFq{p}$.

For $V_\ell$ (resp. $W_\ell$), replace $\kappa_1$ by $\As$
(resp. $\Bs$) in Step 4 as far as reconstruction is concerned.

\medskip
\noindent
{\bf A numerical example:}
Let us give one value for $\ell = 5$. We select $D = -71$ for which
$h(-71) = 7 \geq 5+2$. Consider $p = 1811$. The roots of $H_{-71}(z)$
modulo $p$ are:
$$\mathcal{J}_D = \{313, 1073, 1288, 1312, 1402, 1767, 1808\}.$$
Associated are curves and neighbors for each $j$ value. These can be
found in Table \ref{fig-ell5}. The power sums $\powersums_r$
corresponding to the values are:
$$\begin{array}{r|rrrrrrr}\hline
\EE_i \backslash r & 1 & 2 & 3 & 4 & 5 & 6 \\ \hline
{[1582, 902]} & 0& 105& 1680& 1379& 756& 772 \\
{[1662, 405]} & 0& 527& 1188& 90& 748& 888 \\
{[1451, 1331]} & 0& 1723& 403& 350& 293& 583 \\
{[1013, 747]} & 0& 1133& 18& 1594& 1738& 105 \\
{[224, 753]} & 0& 95& 760& 1790& 1603& 27 \\
{[1128, 1504]} & 0& 155& 669& 1424& 1130& 522 \\
{[91, 725]} & 0& 1793& 1523& 589& 1233& 134 \\
 \hline
\end{array}$$
For instance, $\sigma_6 = c_{6, 0} E_4^3 + c_{6, 1} \Delta$, we need
to solve
$$\left\{
 \begin{array}{ccl}
772 &=& c_{6, 0} 680^3+ c_{6, 1} 1067 \bmod 1811, \\
888 &=& c_{6, 0} 1257^3+ c_{6, 1} 874 \bmod 1811, \\
583 &=& c_{6, 0} 120^3+ c_{6, 1} 363 \bmod 1811, \\
\cdots &\cdots &\cdots \\
 \end{array}
\right.$$
that is $648 E_4^3 + 523 \Delta$. The coefficients are:
\begin{eqnarray*}
\sigma_2 & = & 120 E_4 \\
\sigma_3& = & 960 E_6 \\
\sigma_4& = & 1025 E_4^2 \\
\sigma_5& = & 235 E_6 E_4 \\
\sigma_6& = & 648 E_4^3 + 523 \Delta
\end{eqnarray*}

%%%%%%%%%% SS
\subsubsection{Computing $\mathcal{A}_\ell$ and $\mathcal{B}_\ell$}

In this section, we use $A$ and $B$ instead of $(E_4, E_6, \Delta)$
for ease of presentation.
Once $U_\ell$ is available, we can use equation~(\ref{eqAs0})
in which we plug the series to get
\begin{equation}\label{eqAs}
-3 \ell^4 E_4(q^\ell) \cdot U_\ell'(\sigma_1(q), A(q), B(q)) =
\mathcal{A}_{\ell}(\sigma_1(q), A(q), B(q)).
\end{equation}
Similarly, we would use
$$-2 \ell^6 E_6(q^\ell) \cdot U_\ell'(\sigma_1(q), A(q), B(q)) =
\mathcal{B}_{\ell}(\sigma_1(q), A(q), B(q))$$
to compute $\mathcal{B}_{\ell}$.

We find the coefficients by solving
a linear system (over $\Q$ or using small primes as already
described). We can precompute the powers of the series for $\sigma_1$,
$A$ and $B$ and remark that $U_\ell$ and $\mathcal{A}_\ell$ share a
lot of them. Also, the series $E_4(q^\ell)$ is rather sparse, so that the
product with this quantity is fast.
There is an advantage to compute $\mathcal{A}_\ell$ and
$\mathcal{B}_\ell$ at the same time, sharing as many powers as possible.

Let us turn our attention towards the computations of
$\mathcal{A}_{\ell}$ (resp. $\mathcal{B}_{\ell}$) using
evaluation/interpolation methods. There is nothing special about using
floating point numbers, except that the system we have to solve has size
$O(\ell^2\times \ell^2)$ leading to a $O(\ell^{2 \omega})$ time algorithm.

Some care must be taken when using the isogeny approach. To exemplify
the problem, consider the case $\ell=11$ (similar problems do not
occur for smaller $\ell$'s). The polynomial $\mathcal{A}_{11}$ reads:
$$\mathcal{A}_{11} = a_{1, 1, 0} (X^{11} A) +
+ \cdots + X (a_{11, 6, 0} A^6 + a_{11, 3, 2} A^3 B^2 + a_{11, 0, 4} B^4) +
(a_{12, 5, 1} A^5 B + a_{12, 2, 3} A^2 B^3).$$
To find them, we use a $20 \times 20$ system whose rightmost columns
are
$$M = \left(\begin{array}{cccccc}
\cdots & \kappa_{1, 1} A_1^6 & \kappa_{1, 1} A_1^3 B_1^2 & \kappa_{1,
1} B_1^4 & A_1^5 B_1 & A_1^2 B_1^3 \\
\multicolumn{6}{c}{\cdots} \\
\cdots & \kappa_{1, 12} A_1^6 & \kappa_{1, 12} A_1^3 B_1^2 &
\kappa_{1, 12}
B_1^4 & A_1^5 B_1 & A_1^2 B_1^3 \\
\cdots & \kappa_{1, 13} A_2^6 & \kappa_{1, 13} A_2^3 B_2^2 &
\kappa_{1, 13}
B_2^4 & A_2^5 B_2 & A_2^2 B_2^3 \\
\multicolumn{6}{c}{\cdots} \\
\cdots & \kappa_{1, 20} A_2^6 & \kappa_{1, 20} A_2^3 B_2^2 &
\kappa_{1, 20}
B_2^4 & A_2^5 B_2 & A_2^2 B_2^3 \\
\end{array}\right)$$
The first 12 rows ($12 = \ell+1$) are the $\ell+1$ curves isogenous to
$[A_1, B_1]$. The remaining 8 are taken from the $\ell+1$ curves
isogenous to $[A_2, B_2]$. Consider the product
$$M \times \left(\begin{array}{c}
0 \\ 0 \\ \cdots \\ 0 \\ 1 \\ u \\ v \\ 0 \\ 0
\end{array}\right)
= \left(\begin{array}{c}
\kappa_{1, 1} (A_1^6+u A_1^3 B_1^2 + v B_1^4) \\
{\cdots} \\
\kappa_{1, 12} (A_1^6+u A_1^3 B_1^2 + v B_1^4)\\
\kappa_{1, 13} (A_2^6+u A_2^3 B_2^2 + v B_2^4)\\
{\cdots} \\
\kappa_{1, 20} (A_2^6+u A_2^3 B_2^2 + v B_2^4)\\
\end{array}\right).$$
Given the $A_i$'s and $B_i$'s, we can solve for $u$ and $v$, yielding
a non-zero vector in the kernel of $M$, showing the system is
under-determined.

Fortunately, we can circumvent this problem using the
dual equations
\begin{equation}\label{eqAsdual}
U_\ell'(-\ell \kappa_1, \As, \Bs) (\ell^4 A) =
\mathcal{A}_{\ell}(-\ell \kappa_1, \As, \Bs),
\quad U_\ell'(-\ell \kappa_1, \As, \Bs) (\ell^6 B) =
\mathcal{B}_{\ell}(-\ell \kappa_1, \As, \Bs)
\end{equation}
in Algorithm~\ref{algo-dual}.

\LinesNotNumbered
\begin{algorithm}[hbt]
\caption{Computing $\mathcal{A}_\ell(X, Y, Z) \bmod p$. \label{algo-dual}}
\SetKwProg{Fn}{Function}{}{}
\Fn{\ComputeAMod($\ell$, $D$, $H_{D}(z)$, $p$):}{
\Input{$\ell$ an odd prime, $D$ the discriminant of an imaginary
quadratic order $\mathcal{O}$ of discriminant $D$ with class number
$h(D) \geq \ell+2$; $H_{D}$ is the class
polynomial associated to order $\mathcal{O}$; $p$ prime with
$p\equiv 1\bmod\ell$ and $4 p = t^2 - \ell^2 v^2 D$, $v\not\equiv 0
\bmod \ell$}
\Output{$\mathcal{A}_\ell(X, Y, Z) \bmod p$}

1. $(\EE_i, \{P_{ik}, \EE_{ik}'\}_{1 \leq k \leq
\ell+1})_{1\leq i\leq h} \leftarrow $PartialVolcano($\ell$, $D$,
$H_D(X)$, $p$)\;

2. Inject the values $(\kappa_{1, i}, A_{ik}', B_{ik}')$ in equation
 (\ref{eqAsdual}) to get a linear system with enough rows\;

3. Solve the relevant linear system\;

4. \Return{$\mathcal{A}_\ell$.}
}
\end{algorithm}
Note that the system is $O(\ell^2)\times O(\ell^2)$, leading to a
$O(\ell^{2 \omega})$ time complexity. Adapting this to the case of
$\mathcal{B}_\ell$ is straightforward.

%%%%% S
\section{Computing the isogenous curve {\em \`a la} Atkin}
\label{sct:FrickeA}

The idea is to generalize the approach in
\cite{Atkin88b,Atkin92b,Morain95a}, that is exploit $q$-series
identities to get the parameters $(\kappa, \As, \Bs)$, where we write
$\kappa$ for $\kappa_1$ from now on. As a matter of fact, we need
relations involving $\tilde{E}_{2k} = E_{2k}(q^{\ell})$, from which all
other quantities follows: $\As = -3 \ell^4 \tilde{E}_4$, $\Bs = -2
\ell^6 \tilde{E}_6$.

%%%%%%%%%% SS
\subsection{Properties of $U_\ell$}

We write for readability $U(\kappa, E_4, E_6) = U_\ell(X, E_4, E_6,
\Delta)$ after replacing $\Delta$ by its expression and
$$\partial_\kappa = \frac{\partial U}{\partial \kappa},
\partial_4 = \frac{\partial U}{\partial E_4},
\partial_6 =  \frac{\partial U}{\partial E_6}.$$
and propagate the notation to double derivatives.

The polynomial $U$ is homogeneous with weights, so that
\begin{equation}\label{homU}
(\ell+1) U = \kappa \partial_\kappa + 2 E_4 \partial_4 + 3 E_6
\partial_6.
\end{equation}
Note that partial derivatives of $U$ are also homogeneous polynomials
and we find
\begin{eqnarray}\label{dxx}
\ell \partial_\kappa &=& \kappa \partial_{\kappa\kappa} + 2 E_4
\partial_{\kappa 4} + 3 E_6 \partial_{\kappa 6}, \\
(\ell-1) \partial_4 &=& \kappa \partial_{\kappa 4} + 2 E_4
\partial_{44} + 3 E_6 \partial_{46}, \\
(\ell-2) \partial_6 &=& \kappa \partial_{\kappa 6} + 2 E_4
\partial_{46} + 3 E_6 \partial_{66}.
\end{eqnarray}

%%%%%%%%%% SS
\subsection{Getting the isogenous curve from $U_\ell$}

%%%%%%%%%%%%%%% SSS
\subsubsection{Finding $\tilde{E}_4$}

\begin{proposition}
The value of $\tilde{E}_4$ is given by
$$- \frac{4 \ell (3 E_4^2 \partial_6+2E_6 \partial_4)
- \partial_\kappa (\ell^2 E_4+4  \kappa^2)}
{\ell^4 \partial_\kappa}.$$
\end{proposition}

\medskip
\noindent
{\em Proof:}
We differentiate (using (\ref{operator})) $U(\kappa,
E_4, E_6) = 0$ to get
\begin{equation}\label{eq1}
\kappa' \partial_\kappa + E_4' \partial_4 + E_6' \partial_6 = 0.
\end{equation}
We differentiate (\ref{eqsig}) leading to
$$\kappa' = \frac{\ell}{2} \; (\ell^2 \tilde{E}_2'- E_2') =
\frac{\ell}{24}\; (\ell^2 (\tilde{E}_2^2-\tilde{E}_4) -
(E_2^2-E_4)).$$
Use (\ref{eqsig}) to replace $\ell \tilde{E}_2$ by $2\kappa/\ell+ E_2$
to get 
$$\kappa' = \frac{\ell}{24} \; \left(\frac{4 \kappa^2}{\ell^2} +
\frac{4 \kappa}{\ell} E_2 - (\ell^2 \tilde{E}_4- E_4)\right),$$
that we plug in (\ref{eq1}) together with the expressions for
$\qdiff{E_4}$ and $\qdiff{E_6}$ from equation (\ref{diff46}) to get a
polynomial of degree 1 in $E_2$ whose coefficient of $E_2$ is
$$\kappa \partial_\kappa + 2 E_4 \partial_4 + 3 E_6 \partial_6,$$
which we recognize in (\ref{homU}). Therefore, we get
\begin{equation}\label{eq2}
(\ell+1) U E_2 + \frac{\ell\,\partial_\kappa}{4} \left(\frac{4
\kappa^2}{\ell^2} - (\ell^2 \tilde{E}_4- E_4)\right)  - 2
E_6 \partial_4 - 3 E_4^2 \partial_6 = 0
\end{equation}
from which we deduce $\tilde{E}_4$ since $U(\kappa, E_4, E_6) = 0$. $\Box$

%%%%%%%%%%%%%%% SSS
\subsubsection{Finding $\tilde{E}_6$}

\begin{proposition}
The value of $\tilde{E}_6$ may be written
$$\tilde{E}_6 = -\, \frac{N}{\ell^6 \, \partial_\kappa^3}$$
where $N$ is some polynomial of degree 3 in $\ell$ and given at the
end of the proof.
\end{proposition}

\medskip
\noindent
{\em Proof:}
We differentiate (\ref{eq1}).
\begin{eqnarray}
\label{ea1}
& & \kappa'' \partial_\kappa + {\kappa'}(\kappa' \partial_{\kappa\kappa}
+ E_4' \partial_{\kappa 4} + E_6' \partial_{\kappa 6}) \\ 
\label{ea2} &+&E_4'' \partial_4 + E_4' (\kappa' \partial_{4\kappa} + E_4'
\partial_{44} + E_6' \partial_{46}) \\
 \label{ea3}&+& E_6'' \partial_6 + E_6' (\kappa' \partial_{6\kappa} + E_4'
\partial_{64} + E_6' \partial_{66}) = 0
\end{eqnarray}
We compute in sequence
$$12 E_2'' = 2 E_2 E_2' - E_4' = E_2 (E_2^2-E_4)/6 - (E_2 E_4 - E_6)/3,$$
$$12 {\tilde{E}_2}'' = 2 \tilde{E}_2 \tilde{E}_2' - \tilde{E}_4'
= \tilde{E}_2 (\tilde{E}_2^2-\tilde{E}_4)/6 - (\tilde{E}_2
\tilde{E}_4 - \tilde{E}_6)/3,$$
which give us the value
$$\kappa'' = \frac{\ell}{2} \; (\ell^3 {\tilde{E}_2}''- E_2'')$$
to be used in (\ref{ea1}).
Differentiating relations of (\ref{diff46}), we get
$$E_4'' = \frac{1}{3} \, (E_2' E_4 + E_2 E_4' - E_6'),
\quad E_6'' = \frac{1}{2} \, (E_2' E_6 + E_2 E_6' - 2 E_4 E_4'),$$
to be used in lines (\ref{ea2}) and (\ref{ea3}) respectively. We
replace $\tilde{E}_4$ by its value from (\ref{eq2}), and $\tilde{E}_2$ using
$\kappa = (\ell/2) ( \ell \tilde{E}_2 - E_2)$.
This finally yields an expression
as polynomial in $E_2$: $$C_2 E_2^2 + C_1 E_2 + C_0 = 0.$$
The unknown $\tilde{E}_6$ is to be found in $C_0$ only.

By luck(?)
\begin{proposition}
The coefficients $C_1$ and $C_2$ vanish for a triplet such that
$U_\ell(\kappa, E_4, E_6) = 0$.
\end{proposition}

\medskip
\noindent
{\em Sketch of the proof:} The strategy to
prove this is the same in both cases. Replace
$\partial_{\kappa\kappa}$, $\partial_{44}$ and $\partial_{66}$ by
their values from (\ref{dxx}). Factoring the resulting expressions
yields the same factor $\kappa \partial_\kappa + 2 E_4 \partial_4 + 3
E_6 \partial_6$, which cancels $C_1$ and $C_2$. We add a \verb+SageMath+
script for the convenience of the reader as an appendix to this
work. $\Box$

We are left with 
$$\tilde{E}_6 = -\, \frac{N}{\ell^6 \, \partial_\kappa^3}$$
where $N$ is a polynomial in degree 3 in $\ell$
$$N = -E_6 \partial_\kappa^3 \ell^3 + c_2 \ell^2 + 12
\partial_\kappa^2 \kappa (3 E_4^2 \partial_6+2 E_6 \partial_4)
\ell - \partial_\kappa^3 \kappa^3.$$
The coefficient $c_2$ is heavy looking and we give slightly factored
as a polynomial in $E_4$:
\begin{eqnarray*}
c_2 &=&
18 (\partial_6^2 \partial_{\kappa\kappa}-2 \partial_6 \partial_\kappa \partial_{\kappa 6}+ \partial_{66} \partial_\kappa^2) E_4^4\\
&&+(24 E_6 \partial_4 (\partial_6 \partial_{\kappa\kappa}-
\partial_\kappa \partial_{\kappa 6})
+24 E_6 \partial_\kappa (\partial_{46} \partial_\kappa- \partial_6 \partial_{\kappa 4})+10 \partial_4 \partial_\kappa^2) E_4^2\\
&&+3 \partial_\kappa^2 (7 E_6 \partial_6 - \kappa \partial_\kappa) E_4
+8 E_6^2 (\partial_4^2 \partial_{\kappa\kappa}-2 \partial_4
\partial_\kappa \partial_{\kappa 4}+\partial_{44}
\partial_\kappa^2). \quad \Box
\end{eqnarray*}

%%%%%%%%%%%%%%% SSS
\subsubsection{Numerical example}

Consider $E: Y^2 = X^3+X+3$ over $\GFq{1009}$ and $\ell = 5$. Using
$$U_5(X) = X^6 + 20 X^4 A + 160 X^3 B - 80 X^2 A^2 - 128 X A B - 80
B^2,$$
we select $\kappa = 584$ and compute
$$\partial_\kappa=905, \partial_4=779, \partial_6=140$$
from which $\tilde{E}_4=497$, $\As = 441$. After tedious computations,
we find $\Bs = 997$.

%%%%% S
\section{Implementation and numerical results}
\label{sct:results}

A lot of trials were done using {\sc Maple} programs, some of which
were then rewritten in {\sc Magma} (version 2.26-10), for speed. See
the author's web page. Computing the polynomials for $\ell \leq 100$
takes a few minutes on a classical laptop. Checking them is done using
SEA, as mentioned in~\cite{NoYaYo20}.

We give some examples of the {\em relative height} $\tilde{H}$ for
some of our polynomials. Here $\tilde{H}(P) = H(P)/((\ell+1)\log
\ell)$. Note that these quantities seem to stabilize when $\ell$
increases and are in accordance with Proposition~\ref{sizeUVW}.

$$\begin{array}{|r|r|r|r|r|}\hline
\multicolumn{5}{|c|}{} \\
\ell & \tilde{H}(\Phi_\ell^t) & \tilde{H}(\Phi_\ell^c) &
\tilde{H}(\Phi_\ell^*) & \tilde{H}(U_\ell) \\ \hline
2 & 15.72 & 4.00 & -- & -- \\
3 & 11.14 & 1.51 & -- & 0.32\;\; \\
5 & 11.243 & 0.762 & -- & 0.526  \\ %% 5, s=4
7 & 9.787 & 0.582 & -- & 0.640 \\ %% 7, s=6
11 & 10.130 & 1.842 & 1.120 & 0.670  \\ % 11
13 & 9.565 & 0.367 & 0.941 & 0.688  \\ %% 1
17 & 9.581 & 0.958 & 0.714 & 0.690  \\ %% 5
19 & 9.365 & 0.648 & 0.630 & 0.695  \\ %% 7
23 & 9.438 & 1.995 & 0.419 & 0.698  \\ %% 11
\hline
101 & -- & 1.111 & 0.159 & 0.778 \\ %% 5
103 & -- & 0.740 & 0.249 & 0.779 \\ %% 7
107 & -- & 2.218 & 0.228 & 0.781 \\ %% 11
109 & -- & 0.379 & 0.213 & 0.782 \\ %% 1
\hline
\end{array}$$
Data are computed using the polynomials available in {\sc
Magma}: $\Phi_\ell^c$ is called {\em canonical polynomial} and
$\Phi_\ell^*$ is called {\em Atkin polynomial}. In the case of
$\Phi_\ell^c$, the height depends on $\ell \bmod 12$.
Still, Atkin's minimal
functions remain the best choice for large $\ell$'s.

%%%%% S
\section{Conclusions}

We have given several methods for computing the Fricke and
Charlap-Coley-Robbins polynomials, and manage to adapt known
algorithms for this task. We also included representations as
fractions in polynomials. In some cases, $U_\ell$ has smaller height,
at long as $\ell$ is small.

In isogeny cryptography they are useful
for relatively small $\ell$'s, if we store $(U_\ell, \mathcal{A}_\ell,
\mathcal{B}_\ell)$. If one wants to compute an isogeny, it is enough to
compute a root of $U_\ell$ followed by instantiations of three
polynomials.

Also, we insisted on families of modular forms. Some of
the techniques can be used for {\it ad hoc} forms.

In a follow up work~\cite{Morain24b}, we look at modular polynomials for
cuspidal $\eta$-products, as already described by
Fricke, one of which was recommanded by Atkin to replace the
$U_\ell$ polynomial when $\ell \equiv 11 \bmod 12$.

\bigskip
\noindent
{\bf Acknowledgments.} The author wishes to thank A.~Bostan and
F.~Chyzak for helpful discussions around some aspects of this work; special
thanks to the former for his impressive list of references for the
fast evaluation of hypergeometric functions. Thanks also to L.~De Feo for
his updates on cryptographic applications of isogenies.

\def\noopsort#1{}\ifx\bibfrench\undefined\def\biling#1#2{#1}\else\def\biling#1#2{#2}\fi\def\Inpreparation{\biling{In
  preparation}{en
  pr{\'e}paration}}\def\Preprint{\biling{Preprint}{pr{\'e}version}}\def\Draft{\biling{Draft}{Manuscrit}}\def\Toappear{\biling{To
  appear}{\`A para\^\i tre}}\def\Inpress{\biling{In press}{Sous
  presse}}\def\Seealso{\biling{See also}{Voir
  {\'e}galement}}\def\Editor{\biling{Ed.}{R{\'e}d.}}

\appendix

%%%%%%%%%%%%%%% S
\section{A script to check the computations}

This \verb+SageMath+~\cite{sagemath} script can also be downloaded
from the author's web page.

{\def\baselinestretch{0.9}
\begin{lstlisting}[language=python]
# This script is devoted to the computation and verification of several
# identities related to Fricke polynomials using the notations of the preprint.

# one ring to rule them all
R.<ell,E2,E4,E6,sigma,E4t,E6t,d4,d6,s,ds,ds4,ds6,d46,f,df,df4,df6>
    =PolynomialRing(Rationals(),18)

########## The Fricke case

# returns ell^-4 * ds^-1 * (-12*ell*E4^2*d6 + ...)
def check_E4t():
    E4p=(E2*E4 - E6)/3
    E6p=(E2*E6-E4^2)/2
    E2p=(E2^2-E4)/12
    E2t=(E2+2*sigma/ell)/ell
    sigp=ell/24*(4*sigma^2/ell^2+4*sigma/ell*E2-(ell^2*E4t-E4))
    tmp=sigp*ds+E4p*d4+E6p*d6
    tmp=tmp.numerator()
    print("degree(tmp, E2)=", tmp.degree(E2))
    # check that coeff of E2 is zero
    c1=tmp.coefficient({E2:1})
    # is a multiple of (2*E4*d4 + 3*E6*d6 + f*df), hence 0
    print("c1=", c1.factor())
    # find sigma as a root of constant coefficient
    e4t=tmp.coefficient({E2:0})
    e4t=-e4t.coefficient({E4t:0})/e4t.coefficient({E4t:1})
    # sig contains the value of sigma
    return e4t.factor()

# returns
# ell^-6 * ds^-3 * sigma^-1 * E6^-1 * E4^-1 * (-18*ell^3*E4^5*E6*d6^2*ds+...)
def check_E6t():
    e4t=check_E4t()
    E4p=(E2*E4 - E6)/3
    E6p=(E2*E6-E4^2)/2
    E2p=(E2^2-E4)/12
    E2t=(E2+2*sigma/ell)/ell
    sigp=ell*(4*sigma^2/ell^2+4*sigma/ell*E2-(ell^2*e4t-E4))/24
    # more derivatives
    E4pp=1/3*(E2p*E4+E2*E4p-E6p)
    E6pp=1/2*(E2p*E6+E2*E6p-2*E4*E4p)
    # crucial values
    E4tp=1/3*(E2t*e4t-E6t)
    E2tp=(E2t^2-e4t)/12
    E2pp=1/12*(2*E2*E2p-E4p)
    E2tpp=1/12*(2*E2t*E2tp-E4tp)
    sigpp=ell*(ell^3*E2tpp-E2pp)/2
    # inject diagonal derivatives
    dss = (ell*ds    -2*E4*ds4 -3*E6*ds6)/sigma
    d44 = ((ell-1)*d4-sigma*ds4-3*E6*d46)/(2*E4)
    d66 = ((ell-2)*d6-sigma*ds6-2*E4*d46)/(3*E6)
    # starting point
    tmp=      sigpp*ds+sigp*(sigp*dss+E4p*ds4+E6p*ds6)
    tmp=tmp + E4pp*d4+E4p*(sigp*ds4+E4p*d44+E6p*d46)
    tmp=tmp + E6pp*d6+E6p*(sigp*ds6+E4p*d46+E6p*d66)
    tmp=tmp.numerator()
    c2=tmp.coefficient({E2:2})
    print("E6t.c2=", c2.factor())
    c1=tmp.coefficient({E2:1})
    print("E6t.c1=", c1)
    c0=tmp.coefficient({E2:0})
    e6t=-c0.coefficient({E6t:0})/c0.coefficient({E6t:1})
    return e6t.factor()
\end{lstlisting}}

%%%%%%%%%%%%%%% S
\section{Some values of Fricke/CCR polynomials}
\label{sct:UVWNANB}

%%%%%%%%%% SS
\subsection{Prime indices}

Note there is a sign flip compared to~\cite{NoYaYo20}: They use $B = 3
E_6$, whereas we use $B = -3 E_6$ which is coherent with Atkin's work,
say.

For $\ell = 2$, $U_2(X) = X^3+A X + B$ itself, which is the minimal
polynomial of any of the $2$-torsion points.
For $\ell = 3$, we compute
%% (4, 0, 0) -> 4*4     = 16
%% (3, 1, 0) -> 3*4+1*4 = 16
%% (2, 2, 0) -> 2*4+2*4 = 16
%% (1, 3, 0) -> 1*4+3*4 = 16
%% (1, 0, 2) -> 1*4+2*6 = 16
%% (0, 4, 0) -> 4*4     = 16
%% (0, 1, 2) -> 1*4+2*6 = 16
$$W_3(X, E_4, E_6, \Delta) =
X^4+1464 E_6 X^3+(8760 E_4^3-1185643008 \Delta) X^2+(17504
E_4^3-152195991552 \Delta) E_6 X$$
$$+
11664 E_4^6-1790914074624 \Delta E_4^3-20889728069861376 \Delta^2
;$$
Remark that $U_3 = \psi_3/3$. Also
$$\mathcal{A}_{3}(X, E_4, E_6, \Delta) =
-252 E_4 X^3-720 E_6 X^2-684 E_4^2 X-216 E_4 E_6$$
$$\mathcal{B}_{3}(X, E_4, E_6, \Delta) =
-1464 E_6 X^3-4368 E_4^2 X^2-4344 E_4 E_6 X-1440 E_4^3+746496 \Delta
$$
For $\ell = 5$:
{\small
$$\mathcal{A}_{5} =
-1890 E_4 X^5-18720 E_6 X^4-74160 E_4^2 X^3-146880 E_4 E_6 X^2+(-145440 E_4^3+
199065600 \Delta) X-57600 E_4^2 E_6,
$$
$$\mathcal{B}_{5} =
-31260 E_6 X^5-312480 E_4^2 X^4-1249440 E_4 E_6 X^3+(-2497920
E_4^3+763084800 \Delta) X^2-2496960 E_4^2 E_6 X$$
$$-998400 E_4^4+1725235200 \Delta E_4.
$$
}

$$\begin{array}{|r|r|r|r|r|}\hline
\multicolumn{5}{|c|}{} \\
\ell & \tilde{H}(V_\ell) & \tilde{H}(W_\ell) & \tilde{H}(\mathcal{A}_{\ell})
& \tilde{H}(\mathcal{B}_{\ell}) \\ \hline
5 & 3.266 & 4.336 & 1.063 & 1.268  \\
7 & 3.050 & 4.207 & 0.973 & 1.167  \\
11 & 2.939 & 3.979 & 0.896 & 1.016 \\
13 & 2.856 & 3.969 & 0.864 & 0.983\\
17 & 2.770 & 3.883 & 0.831 & 0.919\\
19 & 2.754 & 3.831 & 0.820 & 0.901\\
23 & 2.723 & 3.764 & 0.799 & 0.869 \\
\hline
101 & 2.471 & 3.527 & 0.801 & 0.818 \\
103 & 2.469 & 3.517 & 0.801 & 0.818 \\
107 & 2.466 & 3.516 & 0.803 & 0.819 \\
109 & 2.467 & 3.515 & 0.803 & 0.819 \\
\hline
\end{array}$$

%%%%%%%%%% SS
\subsection{An example with $N=6$}
\label{ex6}

Let $f = E_4$. We give some details for the computation of the modular
polynomial for $E_4(q^6)$. For each representative $R = [a, b; c, d]$
of a coset, we compute the matrices $[N a, N b; c, d] = U R'$ with $U
\in \Gamma$ and $R' = [A, B; 0, D]$ to which corresponds the conjugate
$f(R' \tau)$. Following Proposition~\ref{propSDf}, we group the
conjugates matrices $R'$ w.r.t. $D$. We find
$$\begin{array}{l|l|l|l}\hline
R & U & R' & f(R' \tau) \\ \hline
{[1, 0; 0, 1]} & {[1, 0; 0, 1]} & {[6, 0; 0, 1]} & 1296+O(z^{36})\\
\hline
{[2, 1; 3, 2]} & {[1, -3; -1, 4]} & {[3, 0; 0, 2]} & 81+19440 z^9+O(z^{18})\\
{[1, 0; 3, 1]} & {[0, 1; -1, 2]} & {[3, 1; 0, 2]} & 81+19440 \zeta_6^3 z^9+O(z^{18}) \\
\hline
{[1, 1; 2, 3]} & {[1, -2; -1, 3]} & {[2, 0; 0, 3]} & 16+3840 z^4+34560 z^8+O(z^{12})\\
{[1, 0; 2, 1]} & {[0, 1; -1, 3]} & {[2, 1; 0, 3]} & 16+3840 \zeta_6^2 z^4+34560 \zeta_6^4 z^8+O(z^{12})\\
{[1, 0; 4, 1]} & {[-1, 2; -2, 3]} & {[2, 2; 0, 3]} & 16+3840 \zeta_6^4 z^4+34560 \zeta_6^8 z^8+O(z^{12})\\
\hline
{[1, 5; 1, 6]} & {[1, -5; -1, 6]} & {[1, 0; 0, 6]} & 1+240 z+2160 z^2+6720 z^3+17520 z^4+O(z^{5})\\
{[1, 0; 1, 1]} & {[0, 1; -1, 6]} & {[1, 1; 0, 6]} & 1+240 \zeta_6 z+2160 \zeta_6^2 z^2+6720 \zeta_6^3 z^3+17520 \zeta_6^4 z^4+O(z^{5})\\
{[1, 1; 1, 2]} & {[0, 1; -1, 6]} & {[1, 2; 0, 6]} & 1+240 \zeta_6^2 z+2160 \zeta_6^4 z^2+6720 \zeta_6^6 z^3+17520 \zeta_6^8 z^4+O(z^5)\\
{[1, 2; 1, 3]} & {[0, 1; -1, 6]} & {[1, 3; 0, 6]} & 1+240 \zeta_6^3 z+2160 \zeta_6^6 z^2+6720 \zeta_6^9 z^3+17520 \zeta_6^{12} z^4+O(z^5)\\
{[3, 1; 5, 2]} & {[-3, 11; -5, 18]} & {[1, 4; 0, 6]} & 1+240 \zeta_6^4 z+2160 \zeta_6^8 z^2+6720 \zeta_6^{12} z^3+17520 \zeta_6^{16} z^4+O(z^5)\\
{[1, 0; 5, 1]} & {[-4, 5; -5, 6]} & {[1, 5; 0, 6]} & 1+240 \zeta_6^5 z+2160 \zeta_6^{10} z^2+6720 \zeta_6^{15} z^3+17520 \zeta_6^{20} z^4+O(z^5)\\
\end{array}$$
The corresponding power sums $S_{D, t}(f)$ have
$q$-expansion in $q^A$ (with $A = N/D$) and rational integer
coefficients. For instance
$$\begin{array}{l|l|l}
D & A & S_1(D) \\ \hline
1 & 6 & 1296+311040 q^6+O(q^{12}) \\
2 & 3 & 162+349920 q^3+2838240 q^6+O(q^9) \\
3 & 2 & 48+322560 q^2+2903040 q^4+8720640 q^6+O(q^8) \\
6 & 1 & 6+362880 q+2943360 q^2+9810720 q^3+O(q^4) \\
\hline
\end{array}$$
Summing all these, we get
$$S_1 = 1512+362880 q+3265920 q^2+10160640 q^3+O(q^4)$$
in which we recognize $1512 E_4$. Finally
{\tiny
$$\Phi[E_4(6 \tau)]
=
X^{12}-1512 E_4 X^{11}+296316 E_4^2 X^{10}+120 (-181381
E_4^3+3782160000 \Delta) X^9$$
$$-270 E_4 (-2610581 E_4^3+45664128000 \Delta) X^8-72 E_4^2 (155634011
E_4^3+97714341360000 \Delta) X^7$$
$$+12 (7370195077 E_4^6+29670256575360000 \Delta
E_4^3+32018727707136000000 \Delta^2) X^6$$
$$-1944 E_4 (170853343 E_4^6-4897879320240000 \Delta
E_4^3+23743887602688000000 \Delta^2) X^5$$
$$+45 E_4^2 (15102174661 E_4^6-9546408010149120000 \Delta
E_4^3+47160528043659264000000 \Delta^2) X^4$$
$$+320 (-2535407921 E_4^9-1030763754097002000 \Delta
E_4^6-347206664136004992000000 \Delta^2
E_4^3+211923078487971840000000000 \Delta^3) X^3$$
$$-373248 E_4 (-1520467 E_4^9-246490368694140000 \Delta
E_4^6-9901962075946860000000 \Delta^2 E_4^3+12428563306452480000000000
\Delta^3) X^2$$
$$+3547348992 E_4^2 (-61 E_4^9-49957886310000 \Delta
E_4^6-15504631732476000000 \Delta^2 E_4^3+15137343021240000000000
\Delta^3) X$$
$$+34828517376 (-E_4^3+54000 \Delta) (-E_4^9+151013228706000 \Delta E_4^6-224179462188000000 \Delta^2 E_4^3+1879994705688000000000 \Delta^3).
$$
}

%%%%% S
\section{Numerical data for the isogeny volcano algorithm}
\label{sct:num}

\begin{table}
$$\begin{array}{|r|c|c|r|}\hline
i & \EE_i = [A_i, B_i]        & \Es_i & \sigma(\Es_i) \\ \hline
1 & [1582, 902] & [594, 422] & 226 \\
&& [1543, 911] & 1542 \\
&& [937, 1244] & 1283 \\
&& [1333, 561] & 1691 \\
&& [879, 342] & 1212 \\
&& [757, 1578] & 1290 \\
\hline
2 & [1662, 405] & [1770, 433] & 529 \\
&& [1439, 1411] & 1536 \\
&& [259, 355] & 1810 \\
&& [382, 1793] & 1733 \\
&& [1472, 543] & 433 \\
&& [413, 1603] & 1203 \\
\hline 
3 & [1451, 1331] & [1096, 1433] & 743 \\
&& [1371, 1367] & 98 \\
&& [1105, 1195] & 207 \\
&& [1657, 1699] & 787 \\
&& [811, 812] & 1769 \\
&& [779, 1311] & 18 \\
\hline
4 & [1013, 747] & [1691, 473] & 1705 \\
&& [509, 342] & 1245 \\
&& [1642, 417] & 1406 \\
&& [127, 765] & 1519 \\
&& [905, 1464] & 145 \\
&& [1277, 254] & 1224 \\
\hline
5 & [224, 753] & [1485, 892] & 1566 \\
&& [823, 1106] & 908 \\
&& [397, 1451] & 1729 \\
&& [131, 673] & 450 \\
&& [654, 1798] & 1353 \\
&& [1805, 1025] & 1238 \\
\hline
6 & [1128, 1504] & [1275, 1672] & 1176 \\
&& [1409, 761] & 1362 \\
&& [907, 1757] & 309 \\
&& [824, 1267] & 781 \\
&& [578, 1320] & 1208 \\
&& [1168, 1207] & 597 \\
\hline
7 & [91, 725] & [1184, 542] & 1284 \\
&& [1753, 297] & 859 \\
&& [1440, 1524] & 1268 \\
&& [421, 410] & 517 \\
&& [1626, 1013] & 245 \\
&& [198, 159] & 1260 \\
\hline
\end{array}$$
\caption{Values for $\ell = 5$ and $p = 1811$. \label{fig-ell5}}
\end{table}

\end{document}